\documentclass[final]{siamltex}

\usepackage{amsmath, amsfonts, amssymb, bm}
\usepackage{graphicx, subcaption}
\usepackage{array}

\allowdisplaybreaks

\def\Dt{\Delta t}
\def\R{{\mathbb R}}
\newcommand{\argmin}{\operatornamewithlimits{argmin}}

\title{Numerical Approach for On-the-Fly Active Flow Control via Flow Map Learning Method}

\author{Xinyu Liu\and Qifan Chen\and Dongbin Xiu\thanks{Department of Mathematics,
		The Ohio State University, Columbus, OH 43210, USA. Emails:
		{\tt liu.12165, chen.11010, xiu.16@osu.edu.} Funding: This 
		work was partially supported by AFOSR FA9550-24-1-0237.}
				}

\begin{document}
\maketitle
\begin{abstract}
We present a data-driven numerical approach for on-the-fly active flow control and demonstrate its effectiveness for drag reduction in two-dimensional incompressible flow past a cylinder. The method is based on flow map learning (FML), a recently developed framework for modeling unknown dynamical systems that is particularly effective for partially observed systems. For active flow control, we construct an FML dynamical model for the quantities of interest (QoIs), namely the drag and lift forces. During offline learning, training data are generated for the responses of drag and lift to the control variable, and a deep neural network (DNN)–based FML model is constructed. The learned FML model enables online optimal flow control without requiring simulations of the flow field. We demonstrate that the FML-based approach can be integrated with existing optimal control strategies, including deep reinforcement learning (DRL) and model predictive control (MPC). Numerical results show that the proposed approach enables on-the-fly flow control and achieves more than $20\%$ drag reduction. By eliminating the need for forward simulations during control optimization, the approach offers the potential for real-time optimal control in other systems.
\end{abstract}


\section{Introduction}

Optimal control of complex systems is a longstanding problem in science and engineering, with applications ranging from fluid mechanics and structural dynamics to chemical processes and energy systems \cite{bryson2018applied, bertsekas2012dynamic, Lions1971}. In many such problems, the governing equations involve nonlinear partial differential equations (PDEs) that are computationally expensive to solve. This creates a fundamental challenge for optimal control, since many control strategies require repeated forward simulations of the underlying PDEs. The difficulty is particularly pronounced when control decisions must be made in real time or across many parameter and control configurations \cite{gunzburger2002perspectives, bewley2001flow}.

Classical approaches to optimal control include adjoint-based methods, model predictive control (MPC), and linear–quadratic formulations; see, for example, \cite{Lions1971, Bewley2001, Biegler2010}. Although these methods are mathematically well established, their practical application to large-scale nonlinear systems often requires intrusive access to the governing equations, accurate adjoint solvers, and repeated solutions of forward and adjoint problems \cite{jameson1988aerodynamic, giles2000introduction}. These requirements limit their applicability in settings involving legacy solvers, complex multiphysics models, or strict latency constraints.

To reduce computational cost, a substantial body of work has focused on reduced-order modeling (ROM) for control. Projection-based techniques such as proper orthogonal decomposition (POD) combined with Galerkin projection have been widely used in flow control and related applications \cite{Holmes1996, Noack2003, Rowley2005}. Extensions incorporating stabilization, closure modeling, and balanced truncation have improved robustness, while data-driven approaches such as dynamic mode decomposition (DMD) and operator inference provide alternative reduced representations \cite{Schmid2010, Peherstorfer2016}. Despite these advances, ROM-based control approaches often face challenges related to model fidelity, stability, and extrapolation beyond the training regime, particularly for strongly nonlinear or transient flows.

In this paper, we take an alternative perspective by formulating the control problem directly at the level of quantities of interest (QoIs), rather than the full system state. In many control problems, objectives, constraints, and feedback signals depend only on a small number of physically meaningful observables, such as lift and drag coefficients, stress extrema, energy norms, or integrated sensor measurements. If the evolution of these QoIs under control can be predicted accurately and efficiently, control synthesis can be decoupled from the full-order state dynamics, significantly reducing computational complexity.

The main contribution of this work is the development of such a framework for efficient optimal control of complex systems. Our approach employs flow map learning (FML), a recently developed data-driven framework for modeling the discrete-time evolution operator of unknown dynamical systems. FML was first introduced in \cite{qin2018data} for modeling autonomous dynamical systems and later extended to parameterized systems \cite{QinChenJakemanXiu_IJUQ}, non-autonomous systems \cite{QinChenJakemanXiu_SISC}, partially observed systems \cite{FuChangXiu_JMLMC20}, PDEs \cite{chen2022deep,WuXiu_modalPDE}, and stochastic dynamical systems \cite{ChenXiu_JCP24, XuEtal_JMLMC24}. A detailed review of FML can be found in \cite{ChurchillXiu_FML2023}. In this work, we adapt FML to construct dynamical models for QoIs as functions of control signals. The resulting model can be integrated with existing optimal control methods without requiring full-order simulations.

As a representative and nontrivial demonstration, we consider the active flow control (AFC) problem governed by the incompressible Navier–Stokes equations \cite{Choi2008, GadElHak1989, cattafesta2011actuators}. Active flow control has long served as a benchmark problem for control methodologies because of its high dimensionality, strong nonlinearity, and sensitivity to actuation \cite{GadElHak2000, Kim2007, Bewley2001, henning2007feedback}. In this study, we focus on drag reduction for two-dimensional flow past a circular cylinder \cite{Rabault_JFM2019, Wang_JFM2023}. The QoIs are the fluid forces acting on the cylinder, represented by the drag and lift forces. Control is applied through small jet streams at the sides of the cylinder and is realized through a single control parameter \cite{williamson1996vortex, zdravkovich1981review, rashidi2016vortex}. By repeatedly executing a high-fidelity CFD solver, we generate training data for the responses of drag and lift to varying control signals. We then construct an FML dynamical model for the evolution of the fluid forces under changing control inputs. Once trained, the FML model allows evaluation of the objective function with respect to the control signal without using the full CFD solver. The resulting non-intrusive surrogate can then be integrated with optimal control methods to realize real-time AFC. In particular, we consider reinforcement learning (RL) and model predictive control (MPC), two fundamentally different control strategies, to demonstrate the flexibility of the FML-based approach.
The overall algorithm enables on-the-fly optimal control: as the CFD solver advances the flow simulation, the FML-based controller provides instantaneous control signals without delay. The CFD solver then updates the boundary conditions accordingly and achieves approximately $20\sim 30\%$ drag reduction.

This paper is organized as follows. Section \ref{sec:setup} describes the problem setup and control objective. Section \ref{sec:prelim} reviews the relevant background, including FML, RL, and MPC. Section \ref{sec:fml} presents the construction of the FML models for the drag and lift subject to excitations, and Section \ref{sec:control} discusses how the learned FML models can be integrated with RL and MPC control strategies. Section \ref{sec:examples} reports the computational results for the drag reduction active flow control problem.
\section{Problem Setup} \label{sec:setup}

The active flow control problem considered in this paper is drag reduction for two-dimensional incompressible flow past a circular cylinder, where the control is realized via two jets of injection and suction at the sides of the cylinder.

\subsection{Computational Setup and Numerical Solver} \label{sec:solver}

\begin{figure}[htbp]
 	\begin{center}
 		\includegraphics[width=0.45\textwidth]{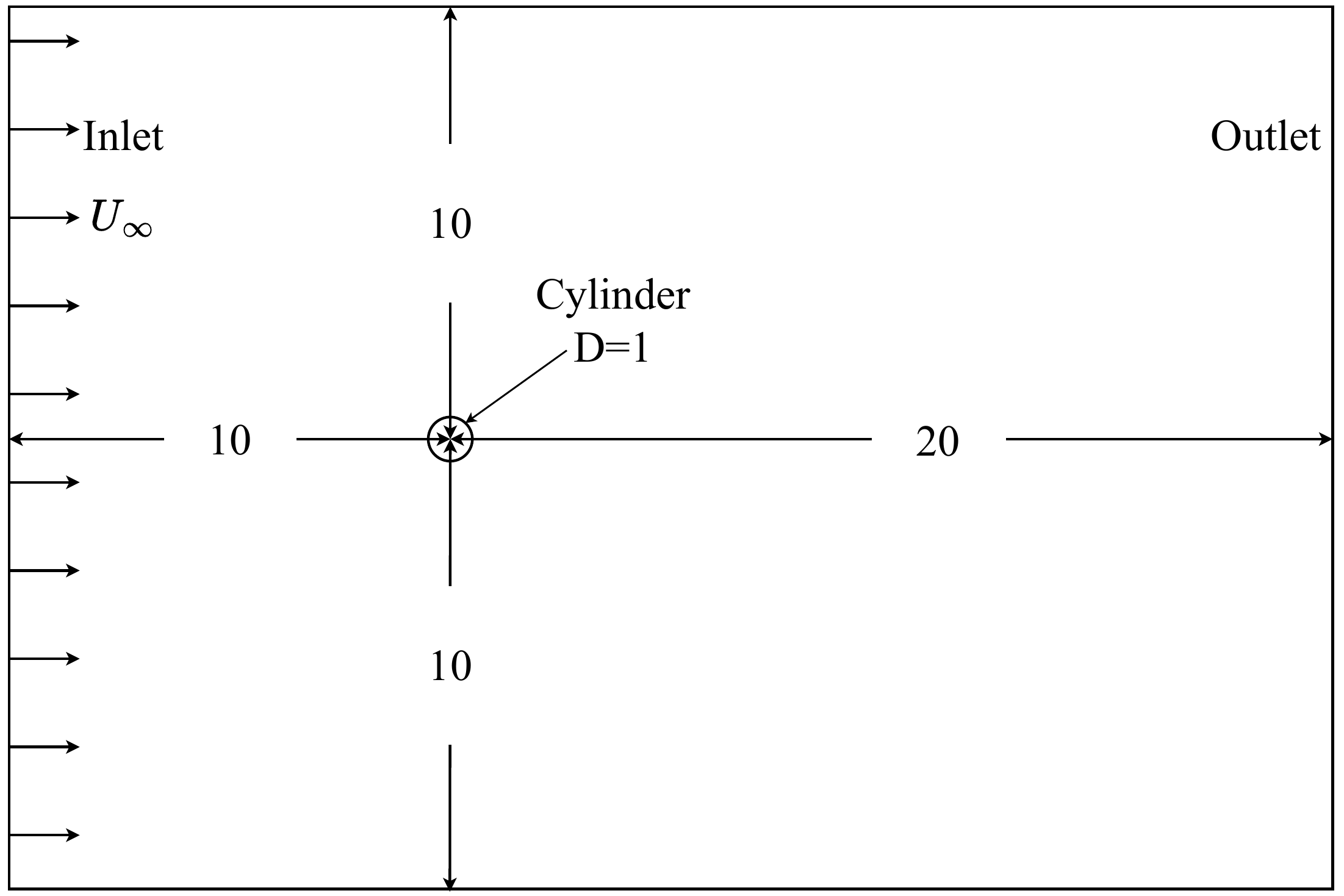}
        \includegraphics[width=0.45\textwidth]{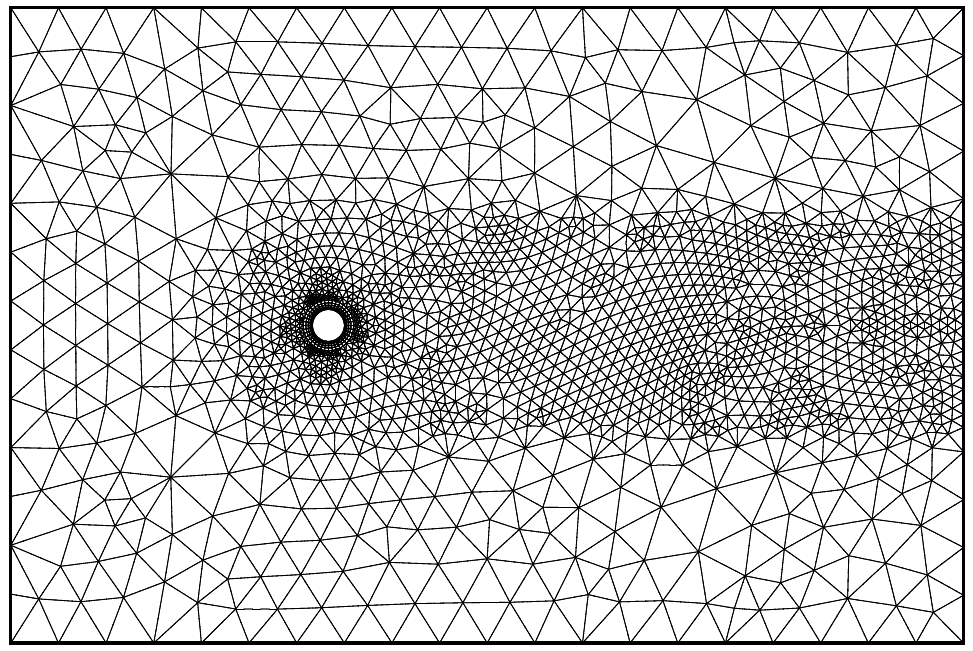}
 		\caption{Computational setup for the flow past cylinder simulation. Left: the computational domain; Right: the computational mesh.}
 		\label{fig:diagram}
 	\end{center}
 \end{figure}

The governing equations are two-dimensional incompressible Navier-Stokes (NS) equations for velocity $\mathbf{u} = (v_x,v_y)^T$ and pressure $p$:
\begin{equation} \label{NS}
    \left\{
      \begin{aligned}
       & \nabla\cdot\mathbf{u} = 0, \\
    &\mathbf{u}_t + (\mathbf{u}\cdot \nabla)\mathbf{u} = -\nabla p +
    \frac{1}{Re}\nabla^2\mathbf{u},
    \end{aligned}
    \right.
\end{equation}
where $Re$ is the Reynolds number
\begin{equation} \label{Re}
Re = \frac{U_\infty L}{\nu},
\end{equation}
with $U_\infty$ being the characteristic velocity, $L$ the characteristic length, and $\nu$ the kinematic viscosity. The computational domain is a rectangular region $(x,y)\in [-10,20]\times [-10,10]$, where a circular cylinder with diameter $D=L=1$ is centered at the origin. For boundary conditions, far-field conditions of $v_x = U_\infty = 1$ and $v_y=0$ are imposed on the left (inlet), top, and bottom boundaries, while an outflow condition of ($\partial_x \mathbf{u}= 0$) is enforced on the right boundary (outlet).
The computational domain is illustrated on the left of Figure~\ref{fig:diagram}.

To solve the problem numerically, we employ Nektar++  (\cite{cantwell2015nektar++,moxey2020nektar++}), an open-source high-order $hp$-spectral element solver. The domain is discretized by an unstructured mesh of $3,613$ elements, shown on the right of Figure~\ref{fig:diagram}. Spectral elements of 6th-order polynomials are used, yielding a total number of $143,880$ degrees of freedom for the discrete system. Second-order implicit–explicit scheme is used for time integration, with a constant step size $\delta t=0.002$. All of these parameter settings were examined numerically to ensure sufficient resolution of the problem.

All simulations are conducted with the ``sudden start" initial state, which sets constant initial velocities $(v_x,v_y)=(1,0)$ inside the domain. We then simulate the flow till it reaches stable and fully developed state. For the Reynolds number considered in this study, $100\leq Re\leq 500$, the flow shall reach well-known quasi-periodic vortex shedding state. Such a state shall be the initial state for our ensuing drag reduction control study.

\subsection{Active Flow Control for Drag Reduction} \label{sec:AFC}

We consider the overall flow force on the cylinder 
\begin{equation*}
    \mathbf{f}(t) = (f_x, f_y)^T(t) = - \oint_S \left(\boldsymbol{\sigma}\cdot\hat{\mathbf{n}}\right) ds,
\end{equation*}
where
$$
    \boldsymbol{\sigma} = -p \mathbf{I} + \nu\left( \nabla \mathbf{u} + \nabla \mathbf{u}^{T} \right)
$$
is the Cauchy stress tensor and the integral is over the cylinder surface $S$ whose outward unit normal vector is 
$\hat{\mathbf{n}}$. In our coordinate system, $f_x$ becomes the drag force and $f_y$ the lift force.

It is a standard practice to normalize the drag and lift forces, resulting in the drag and lift coefficients.
\begin{equation} \label{CD}
    C_D(t) = \frac{f_x}{\frac{1}{2} U_\infty^2 L} = 2f_x(t), 
    \qquad 
    C_L(t) = \frac{f_y}{\frac{1}{2}U_\infty^2 L} = 2 f_y(t),
\end{equation}
where the factor of 2 is the result of our chosen characteristic velocity and length. In this paper, the goal of our active flow control is drag reduction, which shall be implemented by minimizing a cost function involving $C_D$. 
 
\begin{figure}[htbp]
 	\begin{center}
\includegraphics[width=0.45\textwidth]{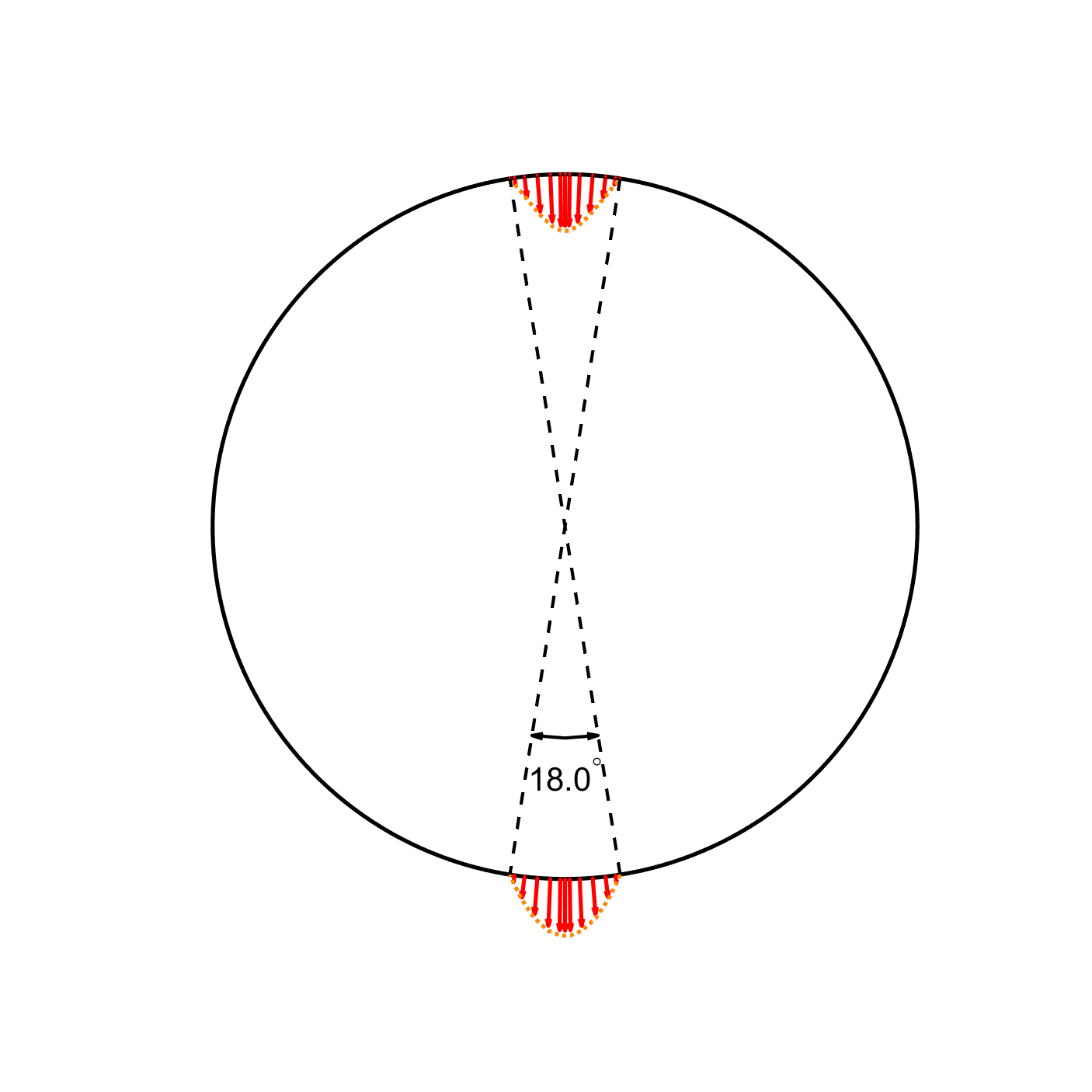}
\caption{Velocity boundary conditions for the two synthetic jets on the cylinder surface.}
 		\label{fig:jets-condition}
        \end{center}
 \end{figure}

To realize drag reduction, we employ an actuator setup similar to that in \cite{Rabault_JFM2019}, where two synthetic jets are placed at the sides of the cylinder as the control. 
More specifically, we introduce an opening of $\omega=18^\circ$ at both the top and bottom of the cylinder boundary, shown in Figure~\ref{fig:jets-condition}.
Two jet flows of identical velocity profile but in the opposite directions are prescribed at the openings, ensuring zero total mass is introduced into the computational domain. 

Specifically, the following velocity profile on the cylinder surface is
\begin{equation}\label{eq:jet_control}
   \mathbf{u}(t) = (v_x,v_y)^T =\frac{2\pi Q(t)}{\omega D^2} \cos \left[ \frac{\pi}{\omega} (\theta - \theta_0) \right]( \cos{\theta}, \sin{\theta})^T, \quad \left|\theta-\theta_0\right|\leq 9^\circ,
\end{equation}
where $\theta$ is the angular coordinate, with $\theta_0=90^\circ$ for the center of the top opening and $\theta_0=270^\circ$ for the center of the bottom opening, $\omega=18^\circ$, and $Q$ is mass flow rate. The optimal control shall be realized via the adjustment of the mass flow rate $Q(t)$. 
        
We define the cost function
\begin{equation} \label{cost}
    J(t) = \langle C_D \rangle_T + \omega_L |\langle C_L\rangle_T|,
\end{equation}
where $\langle\cdot\rangle_T$ stands for moving window time average over window length $T>0$, i.e., 
$$
\langle f \rangle_T(t) = \frac{1}{T}\int_{t-T}^{t} f(\tau) d\tau,
$$
and $\omega_L$ is a parameter. In this paper, we fix $T=5.0$ and $\omega_L=0.2$. These are chosen based on the similar setting from \cite{Rabault_JFM2019}. The goal of the flow control is to achieve active drag reduction via minimization of the cost function,
\begin{equation} \label{control}
    \min_{Q(t)} J(t).
\end{equation}
We remark that there exist a variety of ways for active flow control for flow past a cylinder. The purpose of this paper is to demonstrate the effectiveness of our proposed numerical approach, which is not restricted to the specific choice of the control setup from \cite{Rabault_JFM2019}.

Obviously, the objective function, expressed in terms of the drag and lift coefficients, depends on the underlying flow field, which in turn depends on the $Re$ number. In this paper, we consider two separate types of flow control:
\begin{itemize}
    \item Type-I: The flow field is from a fixed Reynolds number, $Re=300$. This is similar to the work of \cite{Rabault_JFM2019} and most other existing studies.
    \item Type-II: The flow field has an unknown Reynolds number in the range of $100\leq Re\leq 500$. This represents a ``newer" scenario, rarely explored in the literature. In this case, the flow control method needs to implicitly adapt to the flow field whose basic property, manifested via its $Re$ number, is unknown.
\end{itemize}

\section{Preliminaries}\label{sec:prelim}

In this section, we provide preliminary materials essential to our proposed method. These include the fundamentals of FML, as well as the two different optimal control strategies used in this paper, reinforcement learning (RL) and model predictive control (MPC).

\subsection{Flow Map Learning}\label{sec:FML0}

FML is a data driven method for modeling unknown dynamics, first proposed in \cite{qin2018data} for learning autonomous dynamical systems. For an unknown autonomous system, $x_t = f(x)$, $x\in\R^d$, where the governing equation $f:\R^d \to \R^d$ is not known, FML seeks to learn its flow map operator $\Phi:\R^d\to\R^d$ numerically, as opposed to learning $f$. The flow map operator governs the evolution of the solutions, i.e., $x(t_n) = \Phi_{t_n-t_s}(x(t_s))$. Therefore, when data of the solution trajectories are available, we split the trajectory data into pairs of
data separated by one time step $\Delta t$. These pairs of data then satisfy $x(t+\Delta t) = \Phi_{\Delta t}(x(t))$, where $t$ is the ``starting time'' of each data pair. In FML, we then seek a function $N:\R^d\to\R^d$ to minimize the following loss function
$$
\sum_{j=1}^{N_{data}} \left\| x^{(j)}(\Dt) - N(x^{(j)}(0))\right\|^2,
$$
where $N_{data}$ is the total number of such data pairs $\left(x^{(j)}(0), x^{(j)}(\Dt)\right)$, $j=1,\dots,N_{data}$. Note that the ``starting time'' of each pair is set to 0 as in autonomous systems only the relative time matters.
Upon minimizing this loss function, we obtain $N\approx \Phi$ and a predictive FML model 
\begin{equation}\label{fml0}
x_{n+1} = N(x_n),
\end{equation}
which can predict the system behavior when an initial condition $x_0$ is specified.

A notable extension of FML is modeling incomplete system (\cite{FuChangXiu_JMLMC20}). Let $x = (z, w)$, where $z\in\R^m$ are observables and $w\in\R^{d-m}$ missing variables. When trajectory data are available only on $z$, the work of \cite{FuChangXiu_JMLMC20}  established that modeling the dynamics of $z$ can be accomplished without any information of the missing variable $w$. Motivated by the celebrated Mori--Zwanzig formulation (\cite{mori1965, zwanzig1973}), \cite{FuChangXiu_JMLMC20} developed a discrete approximate Mori--Zwanzig formulation and established that a discrete dynamical model of $z$ exists in the following form,
\begin{equation}\label{fml}
      z_{n+1} =
      N(z_n,z_{n-1}, \dots,z_{n-n_M}),  \qquad n\geq n_M,
\end{equation}
where $n_M\geq 0$ is the number of memory terms. In this case, $N:\R^{m\times (n_M+1)}\to\R^m$ is the flow map with memory $n_M$.
The special case of $n_M=0$ corresponds to the learning of a fully observed autonomous system \eqref{fml0}.
The memory-based FML structure \eqref{fml} has been shown to be able to model a wide class of partially observed systems. For a
review of FML modeling of unknown dynamical systems, see \cite{ChurchillXiu_FML2023}. Moreover, a recent work of \cite{ChenXuZhangXiu_JCP26} showed that FML is capable of modeling the dynamics of drag and lift forces for flow past cylinder. This work is directly relevant to the current work.

\subsection{Reinforcement Learning}\label{subsec:RL}
RL is widely applied in control tasks, where an agent learns to optimize its strategy through interactions with an environment \cite{Sutton2018}. Mathematically, this problem is modeled as a Markov Decision Process (MDP). Following the standard definition \cite{Bellman_1957dynamic}, an MDP is characterized by the tuple $(\mathcal{S}, \mathcal{A}, \mathcal{P}, \mathcal{R})$. Here, $\mathcal{S} \subseteq \mathbb{R}^{n_s}$ is the state space, and $\mathcal{A} \subseteq \mathbb{R}^{n_a}$ is the action space. At each step $t_n$, given the current state $s_n \in \mathcal{S}$ and action $a_n \in \mathcal{A}$, the system evolves to next state $s_{n+1}$ according to the transition dynamics $\mathcal{P}: \mathcal{S} \times \mathcal{A} \times \mathcal{S} \to [0, 1]$, which is a conditional density of the next state given the current state and action. Simultaneously, the agent receives a feedback signal from the reward function $\mathcal{R}: \mathcal{S} \times \mathcal{A} \to \mathbb{R}$.

The agent's control strategy is governed by a policy $\pi$, defined as a mapping from the state space to the probability simplex over the action space
\begin{equation}\label{eq:policy}
    \pi: \mathcal{S} \to \Delta(\mathcal{A}).
\end{equation}
Here, $\pi(a_n|s_n)$ denotes the probability of selecting action $a_n$ given state $s_n$. The objective is to find an optimal policy $\pi^*$ that maximizes the expected cumulative discounted reward. This is formulated as
\begin{equation}\label{eq:rl_objective}
    R(\pi) = \mathbb{E}_{\tau \sim (\pi, \mathcal{P})} \left[ \sum_{n=0}^{\infty} \gamma^n \mathcal{R}(s_n, a_n) \right],
\end{equation}
where $\gamma \in [0, 1)$ is the discount factor that determines the present value of future rewards, and $\tau = (s_0, a_0, s_1, a_1, \dots)$ represents a trajectory generated by the policy $\pi$ interacting with the transition dynamics $\mathcal{P}$.

In high-dimensional spaces, solving for the optimal policy $\pi^*$ directly is computationally intractable. Deep RL (DRL) addresses this by approximating the policy function \eqref{eq:policy} using deep neural networks \cite{Mnih_2015Nature}. Consequently, the optimization problem \eqref{eq:rl_objective} transforms into finding the optimal network parameters that maximize $R$.

\subsection{Model Predictive Control}\label{subsec:MPC}
MPC is a well-established framework for controlling constrained dynamical systems \cite{Rawlings_2017MPC}. It determines control inputs by repeatedly solving a finite-horizon optimization problem. Consider a discrete-time system
\begin{equation}\label{eq:mpc_dynamics}
    x_{n+1} = f(x_n, u_n),
\end{equation}
where $x_n \in \mathcal{X} \subseteq \mathbb{R}^{n_x}$ and $u_n \in \mathcal{U} \subseteq \mathbb{R}^{n_u}$ denote the state and control constrained to the admissible sets $\mathcal{X}$ and $\mathcal{U}$. To formulate this optimization problem, a cost function $J$ is defined over a prediction horizon $n_P\geq 1$ as \cite{Mayne_Automatica2000}:
\begin{equation}\label{eq:mpc_cost}
    J(u_n, \dots, u_{n+n_P-1}) = \sum_{k=0}^{n_P-1} \ell(x_{n+k} , u_{n+k}) + V_F(x_{n+n_P}),
\end{equation}
where $\ell(\cdot)$ represents the stage cost and $V_F(\cdot)$ is the terminal cost. At each time step, let $u^*_n$ denote the first element of the optimal sequence that minimizes $J$ subject to the dynamics \eqref{eq:mpc_dynamics} and constraints. The control input is then determined by
\begin{equation}\label{eq:mpc_control}
    u_n = u^*_{n},
\end{equation}
and the process repeats in the next step with the updated state.
\section{FML Modeling for QoI Dynamics under Excitation} \label{sec:fml}

We first discuss how to construct FML models for QoIs dynamics under external excitations. The objective is to construct dynamical model for the drag and lift temporal evolution under arbitrary external input, which is the time dependent boundary condition at the cylinder surface \eqref{eq:jet_control}, without requiring simulation of the flow field. This is an extension of the work of \cite{ChenXuZhangXiu_JCP26}, which demonstrated that FML can model the intrinsic dynamics of drag and lift, without external excitations, for the flow past cylinder problem.

\subsection{FML Models} \label{sec:FML}

We consider the active flow control problem defined in Section \ref{sec:AFC}. Let $V$ denote the QoIs, which are the drag and lift coefficients \eqref{CD}, 
\begin{equation}\label{eq:QoI}
    V(t) = (C_D, C_L)^T(t) \in\mathbb{R}^{n_V}, \qquad n_V=2.
\end{equation}
Let $u$ denote the external control input. In this study, the control is the normalized mass flow rate $Q(t)$ in \eqref{eq:jet_control}, constrained within the interval $[-1,1]$. That is,
\begin{equation}\label{eq:I_u}
    u(t) = Q(t) \in [-1,1]\subset \mathbb{R}^{n_u}, \qquad n_u=1.
\end{equation}

We seek a discrete-time flow map operator $G$ that governs the dynamics of the QoIs $V$ under control excitation $u$. Based on the work of \cite{QinChenJakemanXiu_SISC}, the mapping takes the form
\begin{equation} \label{fml_dim0}
    G:\mathbb{R}^{(n_M+1)\times (n_V+n_u)} \rightarrow \mathbb{R}^{n_V},
\end{equation}
and evolves the QoIs according to
\begin{equation}\label{eq:QoI_dynamic}
    V_{n+1} = G\left(V_n,\cdots,V_{n-n_M}; u_n,\cdots, u_{n-n_M}\right),
\end{equation}
where $n_M\geq0$ is the number of memory steps. With $n_V=2$ and $n_u=1$, we have 
\begin{equation} \label{fml_dim}
    G:\mathbb{R}^{3(n_M+1)} \rightarrow \mathbb{R}^2.
\end{equation}

Note that the Reynolds number $Re$ \eqref{Re} will be treated as a ``hidden" parameter and does not show up explicitly in the model \eqref{eq:QoI_dynamic}. More specifically, for the two types of control problem described in Section \ref{sec:AFC}, we will construct two separate FML models:
\begin{itemize}
    \item FML-1 model: for the Type-I control with flow at the fixed Reynolds number $Re=300$. 
    \item FML-2 model: for the Type-II control where the flow field is from a range of unknown Reynolds numbers $100\leq Re\leq 500$. 
\end{itemize}

In both FML models, the most important parameter is the memory length. To this end, we first determined that a time step $\Delta t=0.1$ is sufficient to resolve the temporal dynamical behavior $V(t)$, the drag and lift coefficients. (The time step used in the Navier-Stokes (NS) solver in Section \eqref{sec:solver}, $\delta t = 0.002$, is suitable for the stability of the NS solver but unnecessarily small for drag and lift dynamics.) We then conduct numerical tests with increasing memory length and determined that the following memory lengths deliver numerically convergent results:
\begin{itemize}
    \item For FML-1, $n_M=20$ which corresponds to a memory length $T_M=2.0$;
    \item For FML-2, $n_M=30$ which corresponds to a memory length $T_M=3.0$. This memory is longer than that of FML-1. It is not surprising because FML-2  ~shall predict the drag and lift for flows with arbitrary and unknown Reynolds number $100\leq Re\leq 500$, as opposed to the fixed $Re=300$ for FML-1.
\end{itemize}

\subsection{Data Generation} \label{sec:sim}

To construct the FML models \eqref{eq:QoI_dynamic}, training data need to be generated by executing the full Navier-Stokes solver. The procedure for generating the training data consists of the following steps:
\begin{itemize}
    \item[(1)] Establish an initial state by simulating the full NS solver under the standard flow past cylinder setting in Section \ref{sec:solver}. For FML-1, the Reynolds number is fixed at $Re=300$; For FML-2, the Reynolds number is randomly sampled in the range of $[100,500]$ with uniform distribution, i.e., $Re\sim \mathcal{U}(100,500)$.
    
    \item[(2)] From the established initial state, continue the NS simulation with randomized time dependent boundary condition at the cylinder surface \eqref{eq:jet_control}. More specifically, the normalized mass flow rate $Q(t)$ in \eqref{eq:jet_control}, which is the control variable $u(t)$ in our setting \eqref{eq:I_u}, is set to be piecewise constant random number following independent uniform distribution in $[-1,1]$ over every $\Delta t=0.1$. The NS simulation is conducted for a time domain $0\leq t\leq T_{sim}$.
    \item[(3)] Record the external excitation signal $u$ and the QoI response over the simulation time $0\leq t\leq T_{sim}$ with the time step $\Delta t=0.1$. This results in a time history data
    \begin{equation} \label{one_data}
       \mathcal{H}\triangleq \{V_{0}, \ldots, V_{N_{step}}; u_{0}, \ldots, u_{N_{step}-1}\}, \quad N_{step}={T_{sim}}/{\Delta t}.
    \end{equation}
    \item[(4)] Repeat the above procedure $N_{sim}>1$ times.
\end{itemize}

Upon finishing the procedure, we obtain a total number of $N_{sim}$ pieces of time history data \eqref{one_data}. These become our raw data set
\begin{equation} \label{data}
    \mathcal{Q} = \left\{\mathcal{H}^{(j)}, j=1,\dots, N_{sim}\right\},
\end{equation}
where $\mathcal{H}^{(j)}$ is the QoI time history data \eqref{one_data} from the $j$-th simulation.

The overall computational cost of obtaining this dataset is the $N_{sim}$ of full NS simulations, each of which lasting $0\leq t\leq T_{sim}$. Although $T_{sim}$ is typically not very long, this procedure is still relatively time consuming. However, we emphasize that this procedure is an offline process. Hereafter, the full NS solver will not be needed anymore.

\subsection{FML Model Construction}

\begin{figure}[htp!]
    \centering
    \includegraphics[width=0.7\textwidth]{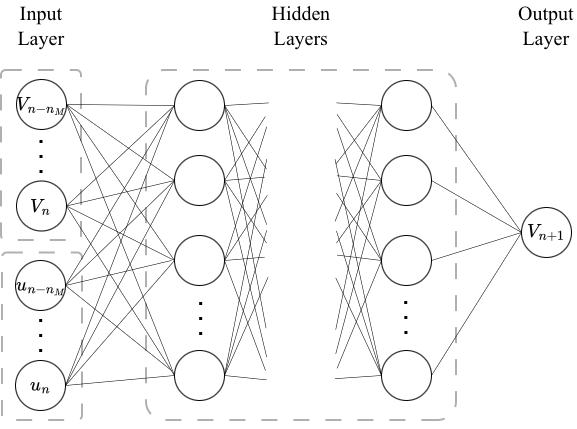}
    \caption{DNN structure for QoI FML model \eqref{eq:QoI_dynamic}.}
    \label{fig:DNN}
\end{figure}

We adopt DNN to express the FML model \eqref{eq:QoI_dynamic}. The DNN structure utilizes a straightforward fully connected feedforward network, as shown in Fig. \ref{fig:DNN}.

To train the DNN, we utilize the dataset \eqref{data}. For each trajectory $\mathcal{H}^{(j)}$, $j=1,\dots, N_{sim}$, in the form of \eqref{one_data}, we 
randomly sample an initial index $n_0$ from $\mathcal{H}^{(j)}$, and then use the segment
$$\left(V_{n_0},\cdots,V_{n_0+n_M}; u_{n_0},\cdots, u_{n_0+n_M-1}\right)
$$
as ``initial condition" to march forward the FML model \eqref{eq:QoI_dynamic} $n_R$ ($n_R\geq 1$) number of times, i.e.,
\begin{equation} \label{fml_train}
\left\{
\begin{aligned}
    &\widetilde{V}_{k+n_M} = G\left(\widetilde{V}_{k+n_M-1},\cdots,\widetilde{V}_{k-1}; u_{k+n_M-1},\cdots, u_{k-1}\right), \quad k=1,\dots,n_R,\\
    &\widetilde{V}_0 = V_{n_0},\quad \cdots,\quad \widetilde{V}_{n_M}=V_{n_0+n_M}, \qquad (\text{Initial conditions})\\
    \end{aligned}
    \right.
\end{equation}
where $k\geq 1$ is the forward marching step. We then 
compute the mean squared loss, averaged over the $n_R$ steps, between the FML prediction $\widetilde{V}$ and the data $V$,
\begin{equation}\label{loss0}
 \frac{1}{n_R} \sum_{k=1}^{n_R} \left\| V_{n_0+n_M+k} - \widetilde{V}_{n_M+k} \right\|^2.
\end{equation}
This is the multi-step loss function we seek to minimize to train the FML DNN models.

The computation of the multi-step loss \eqref{loss0} requires segments of time history data of $V$ \eqref{one_data} of length
\begin{equation}
    n_L = n_M+1+n_R,
\end{equation}
where the first $(n_M+1)$ terms become the inputs of \eqref{fml_train} and the last $n_R$ terms are used to compute the loss \eqref{loss0}. 
Since each $\mathcal{H}^{(j)}$, $j=1,\dots, N_{sim}$ in the dataset \eqref{data} contains $N_{step}+1$ time instances, the random sampling of the initial index $n_0$ thus follows the discrete uniform distribution
\begin{equation}
    n_0 \in \mathcal{U}\left\{0, \dots, N_{step}-n_L\right\}.
\end{equation}

Finally, since the length of the simulated data \eqref{data}, $N_{step}$, is usually (much) larger than $n_L$, we typically conduct $n_B\geq 1$ random sampling of the initial index $n_0$ from each of the $\mathcal{H}^{(j)}$, $j=1,\dots, N_{sim}$, trajectory data in the dataset \eqref{data}. By doing so, the total number of time history segment data used to train the FML models becomes
\begin{equation}
    N_{data} = N_{sim} \cdot n_B.
\end{equation}

\section{FML based Optimal Control} \label{sec:control}

In this section we describe how the FML model \eqref{eq:QoI_dynamic} can be integrated with existing optimal control methods for our active flow control problem. In particular, we discuss its integration with deep reinforcement learning (DRL) and model predictive control (MPC), two widely used and yet fundamentally different control methods, to demonstrate the adaptivity of the FML approach.

\subsection{Deep Reinforcement Learning with FML} \label{sec:RL}

Following the standard notation, our FML enabled DRL here can be expressed via MDP (Markov decision process) tuple $(\mathcal{S}, \mathcal{A}, \mathcal{P}, \mathcal{R})$:
\begin{itemize}
    \item State space $\mathcal{S}$:
The state space $\mathcal{S}$ is defined by the historical trajectory of both QoIs and the applied control signals. Given the memory length $n_M$, the state $s_n$ at each time step $t_n$ is 
\begin{equation}\label{eq:FMLRL_state}
    s_n = \{V_n, \ldots, V_{n-n_M}; u_{n-1},  \ldots, u_{n-n_M}\} \in \mathcal{S}.
\end{equation}

\item Action space $\mathcal{A}$:
The action space $\mathcal{A}$ is the interval $[-1, 1]$, where the control variable $u$ \eqref{eq:I_u} is chosen via a policy $\pi_\theta:\mathbb{R}^{3n_M+2}\to\mathbb{R}$. We also apply action smoothing for smooth control transitions. Specifically, let 
\begin{equation} \label{pi}
    a_n \sim \pi_\theta(\cdot|s_n)
\end{equation}
be a raw action at time $t_n$ generated by the policy $\pi_\theta$, the control at $t_n$ is then defined as
\begin{equation}\label{eq:FMLRL_control}
    u_n = (1 - \alpha)u_{n-1} + \alpha a_n, \qquad \alpha=0.5.
\end{equation}

\item State transition dynamics $\mathcal{P}$: This is the dynamical model of the QoI, which is the FML model \eqref{eq:QoI_dynamic} and constructed via the procedure described in the previous section. Given the state \eqref{eq:FMLRL_state} and the new control \eqref{eq:FMLRL_control}, the FML model \eqref{eq:QoI_dynamic} is fully determined and can be used to evaluate $V_{n+1}$ and obtain the new state $s_{n+1}$.

\item Reward $\mathcal{R}$: The reward function is the negative of the cost function \eqref{cost}, as DRL seeks to maximize the reward. Specifically, at each time $t_n$, the instantaneous reward is
\begin{equation}\label{eq:FMLRL_reward}
    \mathcal{R}_n = -J(t_n) = -\langle C_D\rangle_T(t_n) - \omega_L \left|\langle C_L \rangle_T(t_n)\right|.
\end{equation}

\end{itemize}

 The goal is to learn an optimal policy $\pi_\theta:\mathbb{R}^{2n_M+1}\to\mathbb{R}$ to maximize the expected cumulative reward \eqref{eq:rl_objective}.
  In this paper, we express the policy $\pi_\theta$ via a DNN with $(3n_M+2)$ input nodes, which is the dimension of $s_n$ in \eqref{eq:FMLRL_state}, 1 output node, 2 hidden layers of width 512, and hyper-parameters $\theta$ which are initialized randomly.
  To compute the expected cumulative reward \eqref{eq:rl_objective}, we again start the flow at stable initial state without control, as described in Step (1) in Section \ref{sec:sim}, march forward the learned FML model \eqref{eq:QoI_dynamic} to a finite steps $n_{RL}=200$ (equivalent of a forward marching time $T_{RL}=20.0$), where each step requires a sampling of the policy $\pi_\theta$ \eqref{pi}. The instantaneous reward \eqref{eq:FMLRL_reward} is evaluated at each time step to compute the expected cumulative reward
  \begin{equation} \label{reward}
      R(\pi_\theta) = \mathbb{E}\left[ \sum_{n=0}^{n_{RL}} \gamma^n \mathcal{R}_n \right], \qquad \gamma = 0.99, \quad n_{RL}=200,
  \end{equation}
  where the expectation is evaluated by repeating the above procedure
$5,000$ times such that there are $10^6$ random samples of $\pi_\theta$.
The optimization problem is solved by the proximal policy optimization (PPO) algorithm 
\cite{Schulman_2017PPO}.

The use of the FML model \eqref{eq:QoI_dynamic} renders the entire DRL learning process free of the full NS simulation. As a result, the expected cumulative reward \eqref{reward} is computed via a large number of samples, $10^6$ in this case. In the traditional PDE-based optimal control, this would require $10^6$ PDE simulations, which could be difficult, if not impossible, to accomplish.

\subsection{Model Predictive Control with FML} \label{sec:MPC}

MPC is inherently an online control method, i.e., there is no offline learning. To integrate the learned FML QoI model \eqref{eq:QoI_dynamic} into the flow control, we again start the flow at a full developed stable initial state, as described in Step (1) in Section \ref{sec:sim}. Then, at any time $t_n$ with state $s_n$ \eqref{eq:FMLRL_state}, we conduct forward FML march for a fixed $n_{P}\geq 1$ prediction horizon. In this paper, we fixed $n_{P}=20$, which is equivalent to a time horizon of $T_{P}=2.0$. 
That is, for any fixed $t_n$, given the current state $s_n$ \eqref{eq:FMLRL_state}, we compute
\begin{equation} \label{mpc_fml}
V_{n+k+1} = G\left(V_{n+k},\cdots,V_{n-n_M+k}; u_{n+k},\cdots, u_{n-n_M+k}\right), \quad k=0,\dots,n_{P}-1.
\end{equation}
The cost function \eqref{cost} is then evaluated at the last step
\begin{equation} \label{mpc_obj}
    J(t_n+T_P) = \langle C_D \rangle_T(t_n+T_P)+ \omega_L |\langle C_L \rangle_T(t_n+T_P)|.
\end{equation}
We then seek to minimize this cost function
\begin{equation} \label{mpc_opt}
    \left( u_n^*,\cdots, u_{n+n_P-1}^* \right) = \argmin_{u_n,\cdots, u_{n+n_P-1}}  J(t_n+T_P),
\end{equation}
and set the control at the immediate next step to be
\begin{equation} \label{mpc_un}
    u_n = u_n^*.
\end{equation}
The procedure is then moved to the next time $t_{n+1}$ and continues till the end.

In the traditional MPC, the optimization step \eqref{mpc_opt}, which is carried out by an iterative algorithm such as gradient descent, requires evaluations of the QoIs repetitively. Each evaluation of the QoIs is a full-scale NS solver in this case. The use of the FML model \eqref{mpc_fml} completely avoids the full scale NS solver. Since the FML model \eqref{mpc_fml} is an ultra small 2-dimensional discrete dynamical system, it can be evaluated with negligible cost and seamlessly embedded into the optimization step \eqref{mpc_opt}.

\section{Active Flow Control for Drag Reduction} \label{sec:examples}

In this section, we present the implementation and computational results of the active drag reduction problem and demonstrate the feasibility of on-the-fly control by using  the proposed FML model for the QoIs \eqref{eq:QoI_dynamic}.

\subsection{Numerical Procedure}

Our FML assisted optimal control is implemented in ``real time", in the sense that the feedback active flow control is conducted on-the-fly as the flow field evolves. Specifically, the implementation procedure is as follows:
\begin{itemize}
    \item The full NS solver starts with a fully developed stable state, as described in Step (1) of Section \ref{sec:sim}. At this stage the corresponding control signals are null signals (i.e., uncontrolled).
    \item During the active control stage, the flow field evolves by following the NS simulation with the time dependent boundary condition at the surface cylinder \eqref{eq:jet_control}. The QoIs (drag and lift coefficients) \eqref{eq:QoI} are continuously evaluated by the NS solver. This produces time history data for the QoIs and the corresponding control signals and defines the current state $s_n$ \eqref{eq:FMLRL_state}.
    \item The current state $s_n$ is then fed into the optimal controller. As described in the previous section, we investigate two separate controllers: DRL and MPC.
    \begin{itemize}
        \item In DRL, a raw action is drawn from the learned optimal policy $\pi_\theta$ \eqref{pi} and gives the control signal $u_n$ \eqref{eq:FMLRL_control}.
        \item In MPC, the control signal is computed by \eqref{mpc_un}.
    \end{itemize}
    \item The optimal control value $u_n$ updates the boundary condition at the cylinder surface \eqref{eq:jet_control}. 
    \item The full NS continues with the updated boundary condition.
\end{itemize}

It is evident that while the full NS simulation is conducted to evolve the flow field, the computation of the control signal, via either DRL or MPC, requires only evaluations of the FML model for the QoIs \eqref{eq:QoI_dynamic} and is free of any NS simulations. The cost of computing of the optimal control signal is so minimal that it incurs negligible delay in the full NS simulations. Thus, our active flow control via FML is conducted on-the-fly and in sync with the evolution of the flow field.

\subsection{Open Loop Accuracy Validation}

We first validate the accuracy of the learned FML model \eqref{eq:QoI_dynamic}, which includes two separate FML models for the two types of flow control problems described in Section \ref{sec:AFC}. A summary of the two FML models \eqref{eq:QoI_dynamic} is as follows.
\begin{itemize}
    \item FML-1: Type-I control for flow field at fixed $Re=300$. Here the memory term is $n_M=20$, which defines a mapping $\mathbb{R}^{63}\to\mathbb{R}^2$ \eqref{fml_dim}. We use 4 hidden layers with equal width of $50$ nodes. Therefore, the FML-1 has layer widths $\{63,50,50,50,50,2\}$. We will primarily use FML-1 for DRL control (Section \ref{sec:RL}), where its accuracy needs to be maintained for the forward marching time $T_{RL}=20.0$.
    \item FML-2: Type-II control for flow field at unknown $100\leq Re\leq 500$. Here the memory term is $n_M=30$, which defines a mapping $\mathbb{R}^{93}\to\mathbb{R}^2$ \eqref{fml_dim}. We use 4 hidden layers with equal width of $80$ nodes. Therefore, the FML-2 has layer widths $\{93,80,80,80,80,2\}$. We will primarily use FML-2 for MPC control (Section \ref{sec:MPC}), where its accuracy needs to be maintained for the prediction time horizon $T_{P}=2.0$.
\end{itemize}

Both models were trained using the Adam optimizer to minimize the multi-step loss function \eqref{loss0} over $n_R=3$ recurrent steps. Training was conducted for 100,000 epochs with batch size 4,096. An exponential learning rate schedule was used with initial rate $0.0005$ and decay factor $0.9999$. 

\begin{figure}[ht!]
    \centering
    {\includegraphics[width=0.475\textwidth]{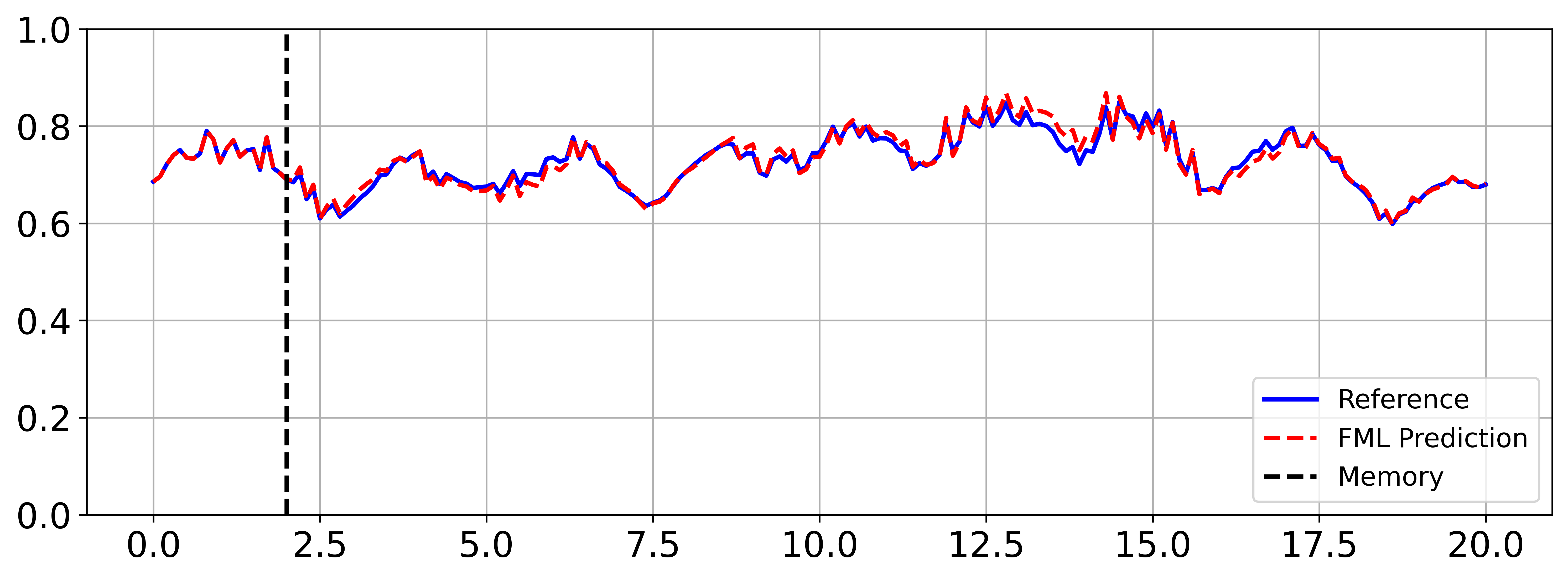}
	\includegraphics[width=0.475\textwidth]{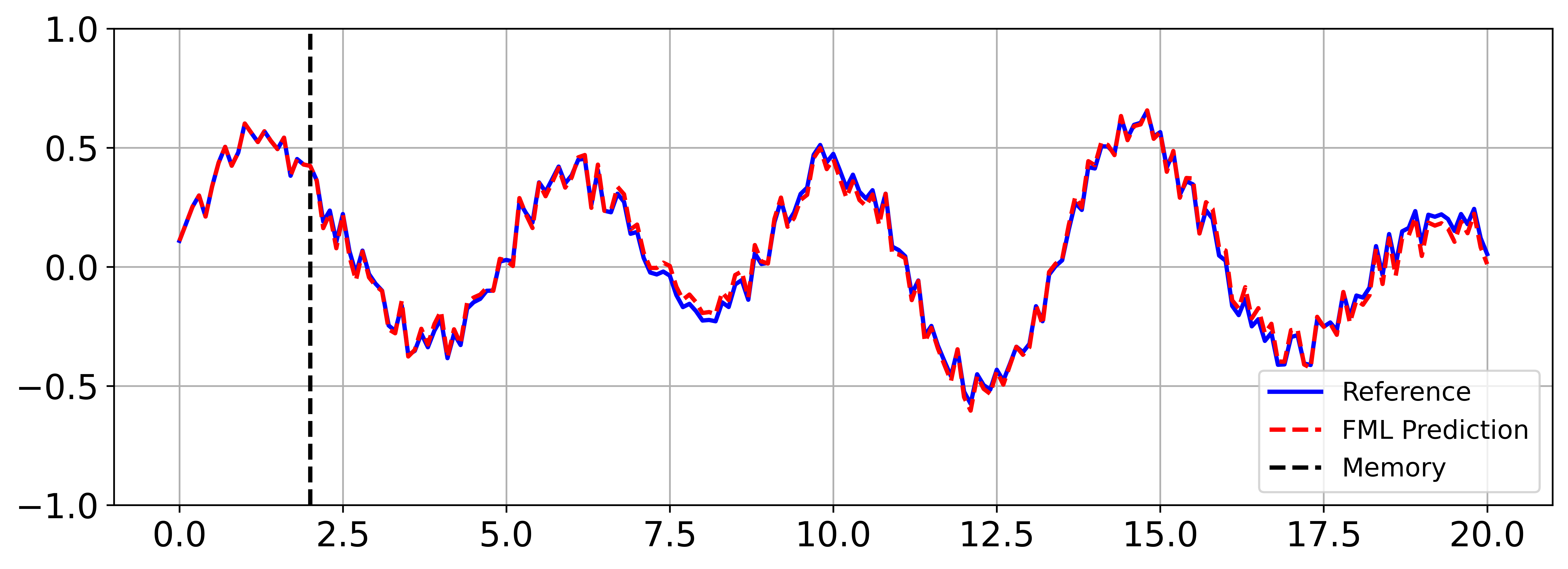}}
    \caption{Open-loop validation for the trained FML-1 model at $Re=300$. Left: drag coefficient ($C_D$); Right: lift coefficient ($C_L$).}
    \label{fig:FML1_open_loop_val}
\end{figure}

To validate the predictive accuracy of the trained FML models, we performed open-loop tests, where the time dependent boundary conditions at the cylinder surface \eqref{eq:jet_control} are subject to randomly generated external signals. (No feedback control is enforced.) 
In Fig. \ref{fig:FML1_open_loop_val}, we show the predictions of the drag and lift coefficients by the FML-1 model, along with the reference solution which was obtained by running the full scale NS solver subject to the same randomized boundary conditions. We observe excellent agreement between the FML predictions and the reference solutions for up to $T_{RL}=20$, which is required by the Type-I control using DRL.

\begin{figure}[ht!]
    \centering
    {\includegraphics[width=0.475\textwidth]{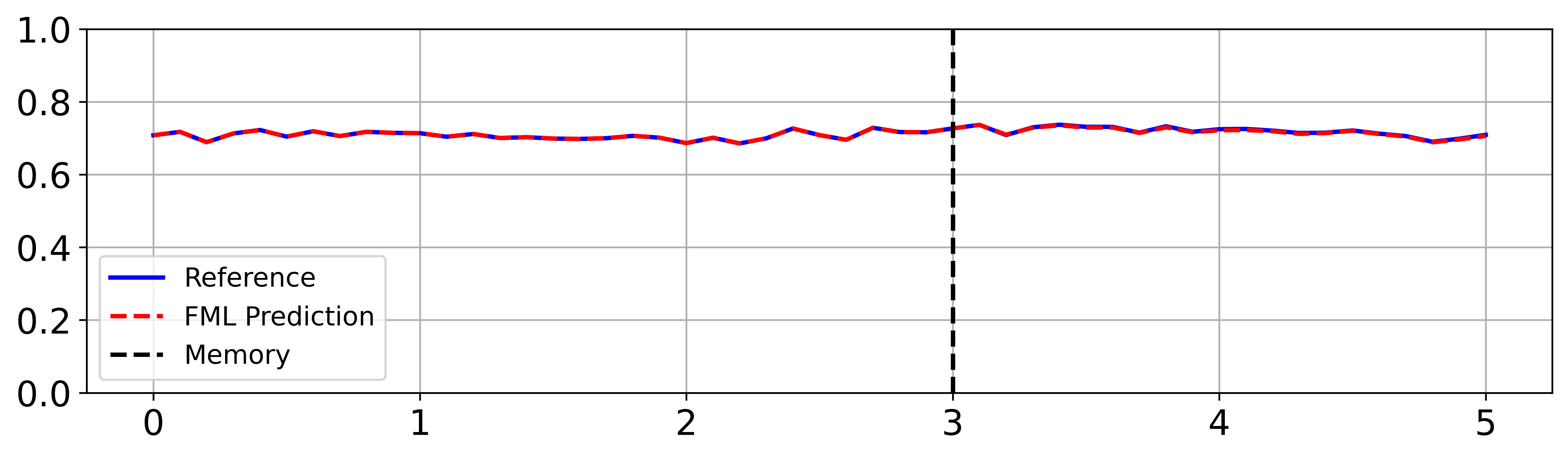}
	\includegraphics[width=0.475\textwidth]{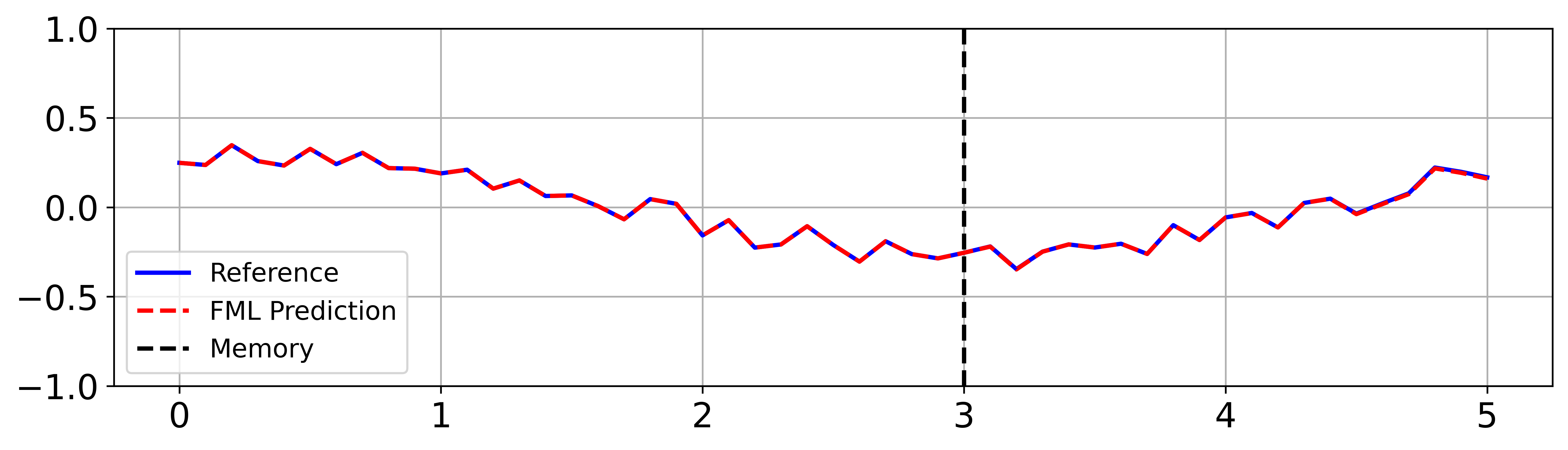}}
    {\includegraphics[width=0.475\textwidth]{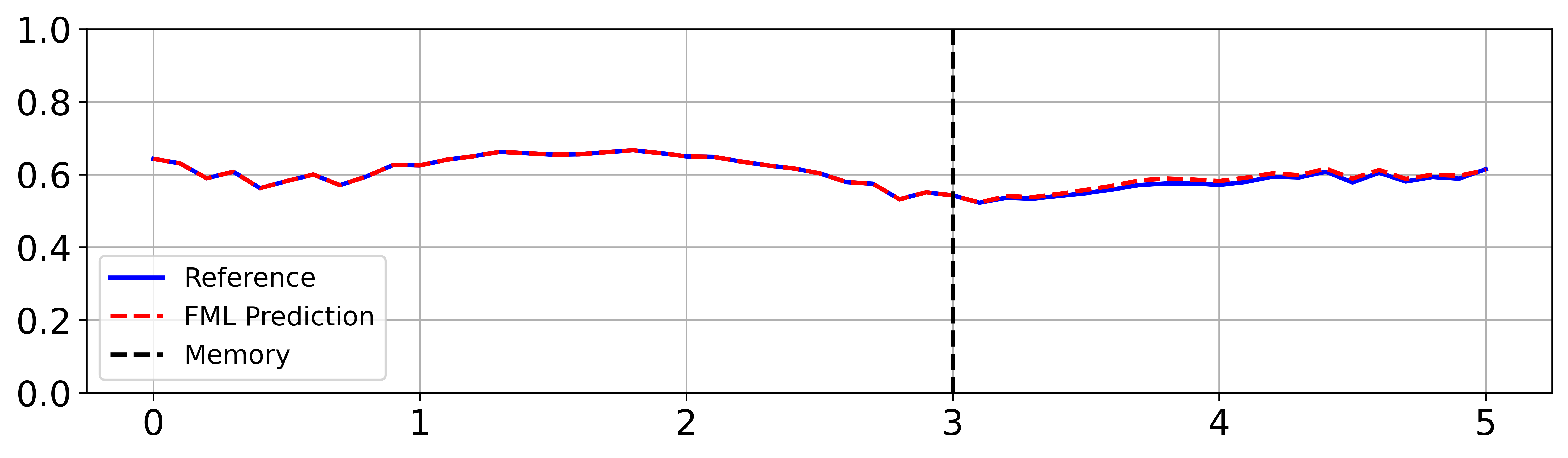}
	\includegraphics[width=0.475\textwidth]{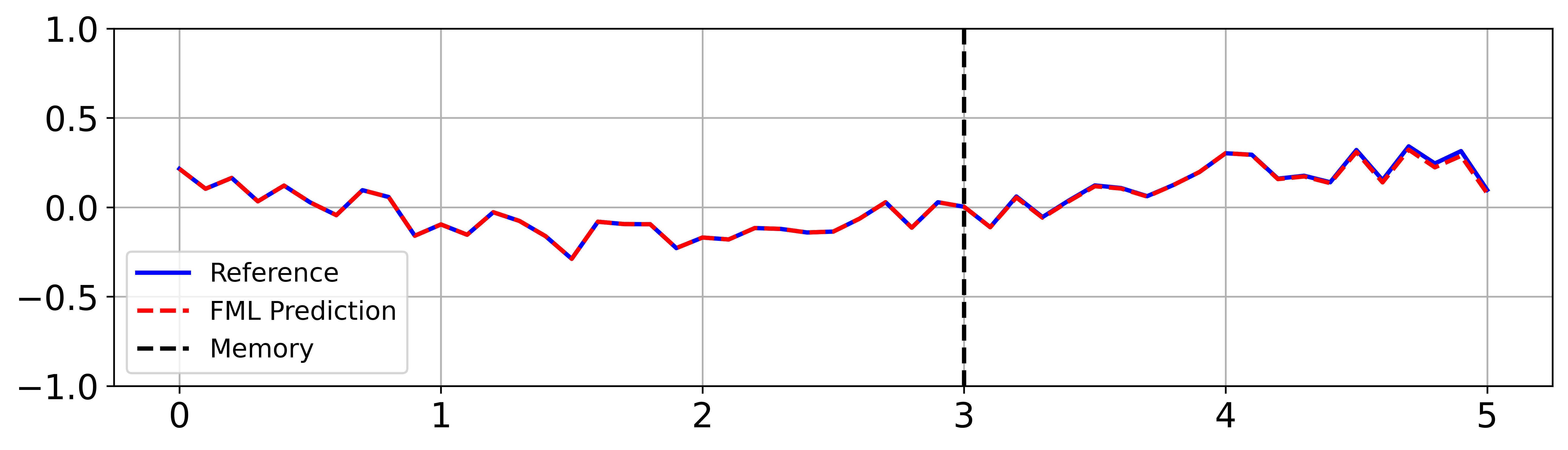}}
    {\includegraphics[width=0.475\textwidth]{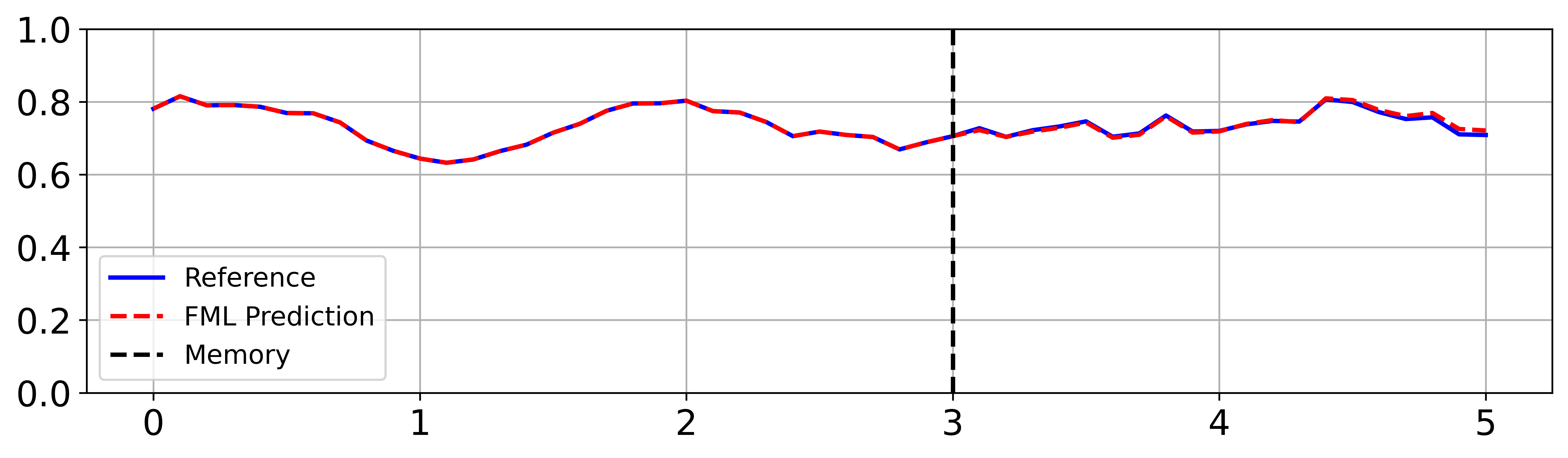}
	\includegraphics[width=0.475\textwidth]{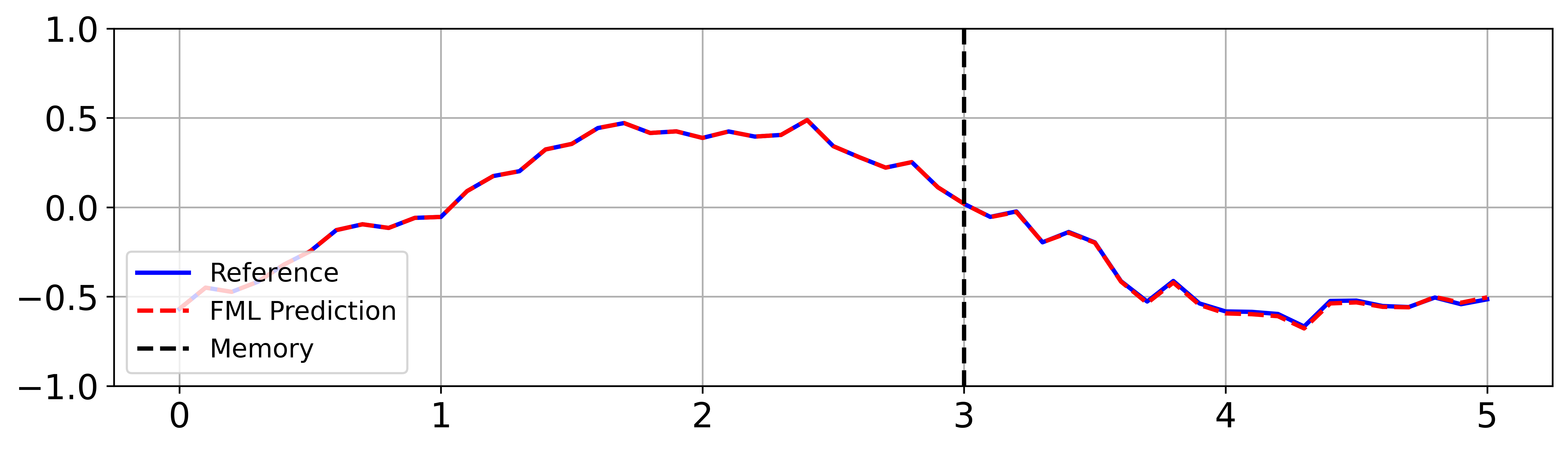}}
    \caption{Open-loop validation of the trained FML-2 model. Left column: drag coefficient ($C_D$); Right column: lift coefficient ($C_L$). Top row:  $Re = 118.62$; Middle row: $Re = 303.10$; Bottom row: $Re = 444.22$. The $Re$ numbers are randomly generated from $Re \in [100,500]$ and unknown to the FML-2 model.}
    \label{fig:FML2_open_loop_val}
\end{figure}

In Fig. \ref{fig:FML2_open_loop_val}, we show the predictions by the FML-2 model, along with the reference solution obtained by the full-scale NS solver. The three rows of results are for, from the top to bottom, $Re=118.62$, $Re=303.10$, and $Re=444.22$. These are randomly generated Reynolds numbers $100\leq Re\leq 500$ and are ``hidden" from the FML-2 model. The only inputs to the FML-2 model are the ``initial conditions", which are the drag and lift coefficient history curves to the left of the vertical dashed lines. The FML-2 model is able to produce accurate predictions of the drag and lift coefficients, compared to the reference solutions obtained by the full-scale NS solver, for the prediction horizon of $T_P=2.0$, which is required by the Type-II optimal control using MPC.

We remark that with sufficient training data and proper training, FML models are capable of producing accurate predictions for (much) longer time horizon. For example, the work of \cite{ChenXuZhangXiu_JCP26} showed that the FML models were accurate for up to $T=200.0$ for predictions of drag and lift for the (uncontrolled) flow past cylinder problem.  Here, we restrain our accuracy tests to time horizons of $T_{RL}=20.0$ for FML-1 and $T_P=2.0$ for FML-2, because they are sufficient for the Type-I control with DRL and the Type-II control with MPC, respectively.

\subsection{Type-I Active Flow Control}
This is for flow at fixed $Re=300$, where we use the trained FML-1 model (Section \ref{sec:FML}) along with DRL-PPO control (Section \ref{sec:RL}). For comparison purpose, we also conducted active control with MPC with the same FML-1 model.

\begin{figure}[ht!]
    \centering
    \includegraphics[width=\textwidth]{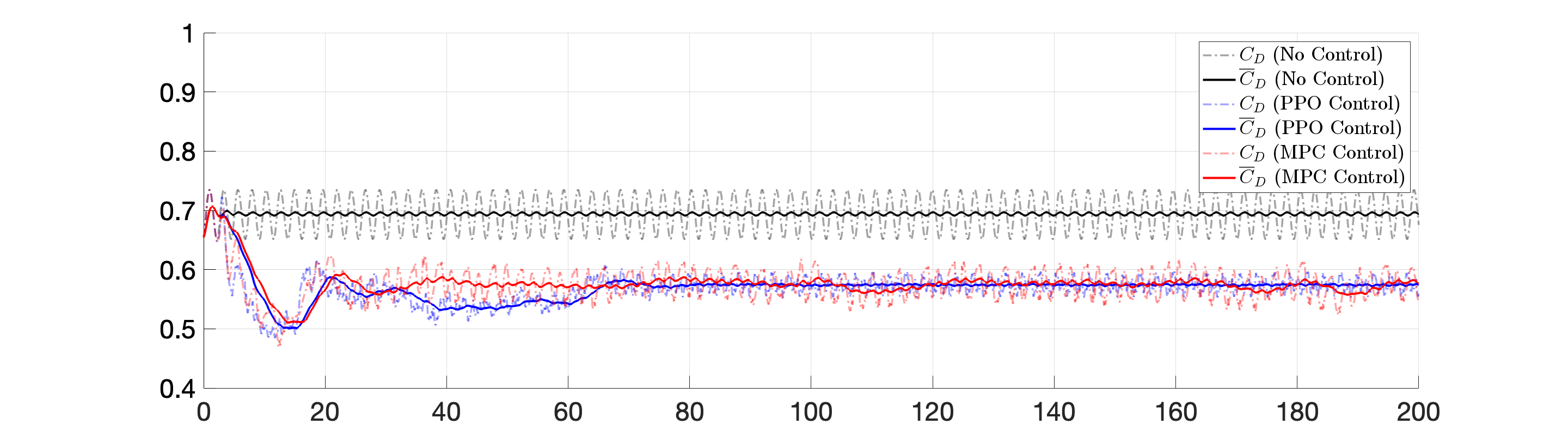}
    \includegraphics[width=\textwidth]{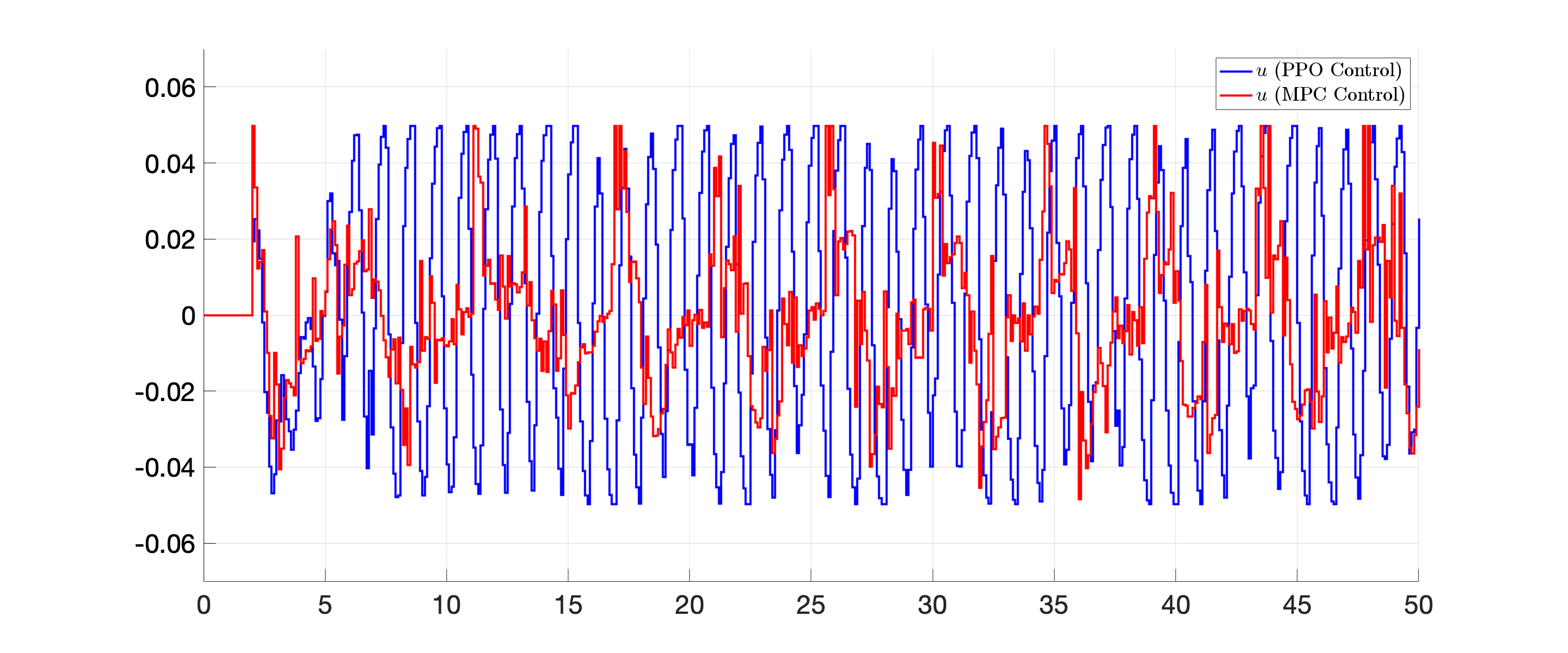}
    \includegraphics[width=\textwidth]{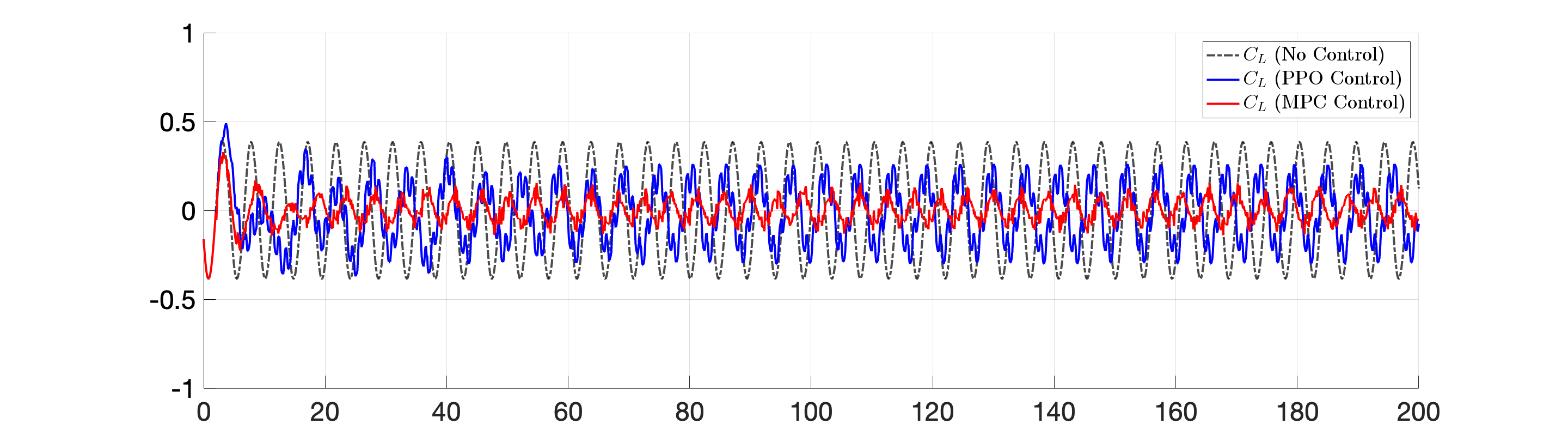}
    \caption{Type-I optimal control at $Re=300$. Top: drag coefficient; Middle: control signal; Bottom: lift coefficient. Dashed lines are instantaneous values; Solid lines are averaged values over moving window of size $5.0$.}
    \label{fig:closed_loop_case1}
\end{figure}

\begin{figure}[ht!]
    \centering
    
    \begin{tabular}{ m{0.42\textwidth} @{}c m{0.42\textwidth} @{}c@{} }
        
        \includegraphics[trim=8cm 2cm 8cm 2cm,clip,width=0.42\textwidth]{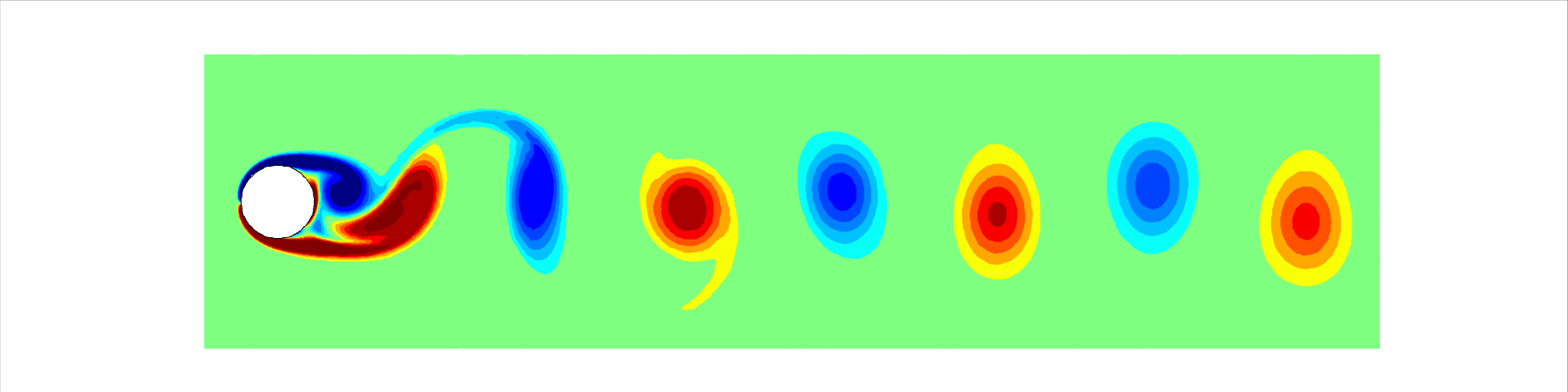} &
        \rotatebox[origin=c]{90}{$t=21.8$}&        
        \includegraphics[trim=8cm 2cm 8cm 2cm,clip,width=0.42\textwidth]{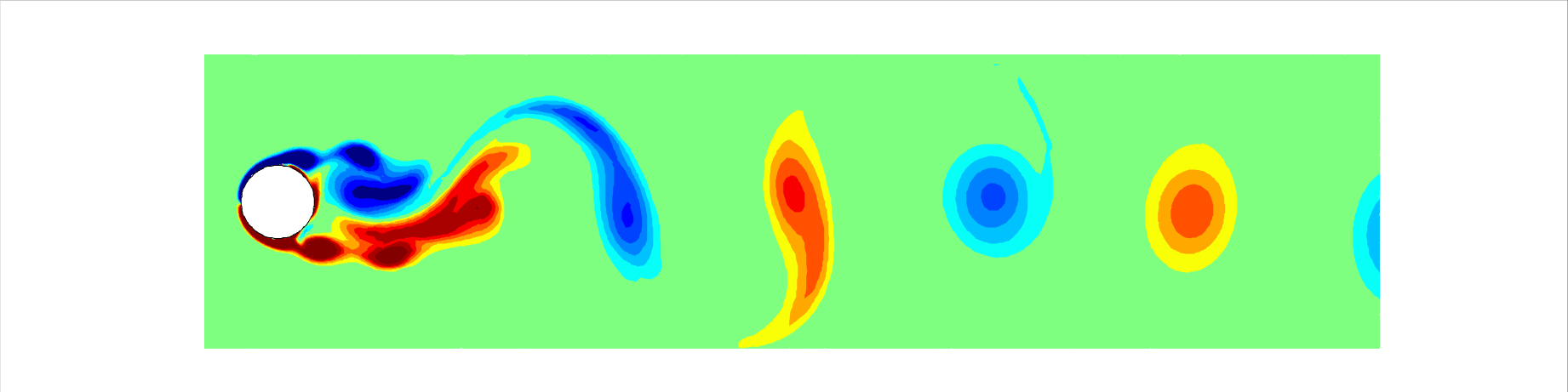}         &
        \rotatebox[origin=c]{90}{$t=80.5$}
        \\[1em]
        
  \includegraphics[trim=8cm 2cm 8cm 2cm,clip,width=0.42\textwidth]{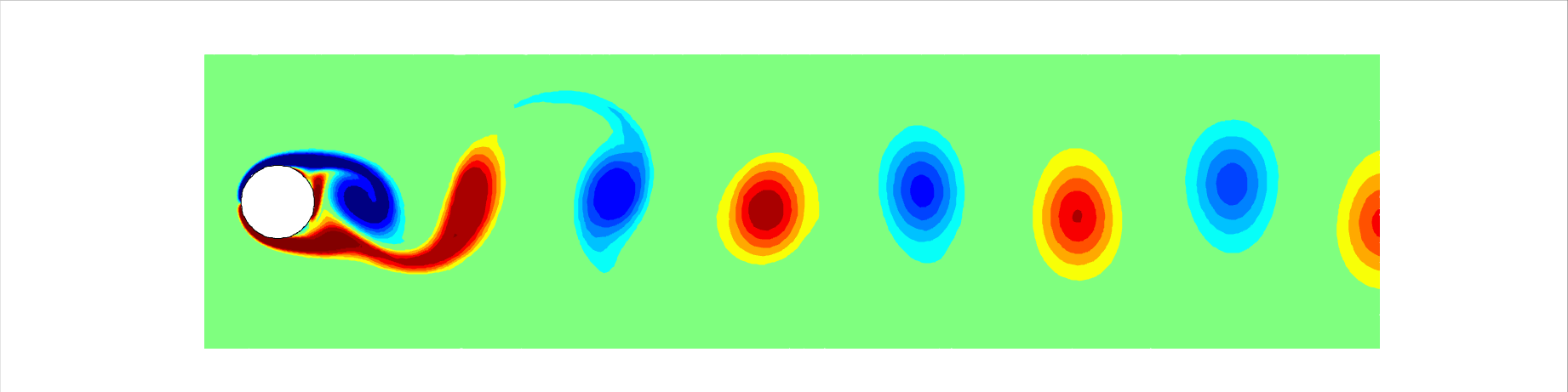} &
        \rotatebox[origin=c]{90}{$t=23.0$}&        
        \includegraphics[trim=8cm 2cm 8cm 2cm,clip,width=0.42\textwidth]{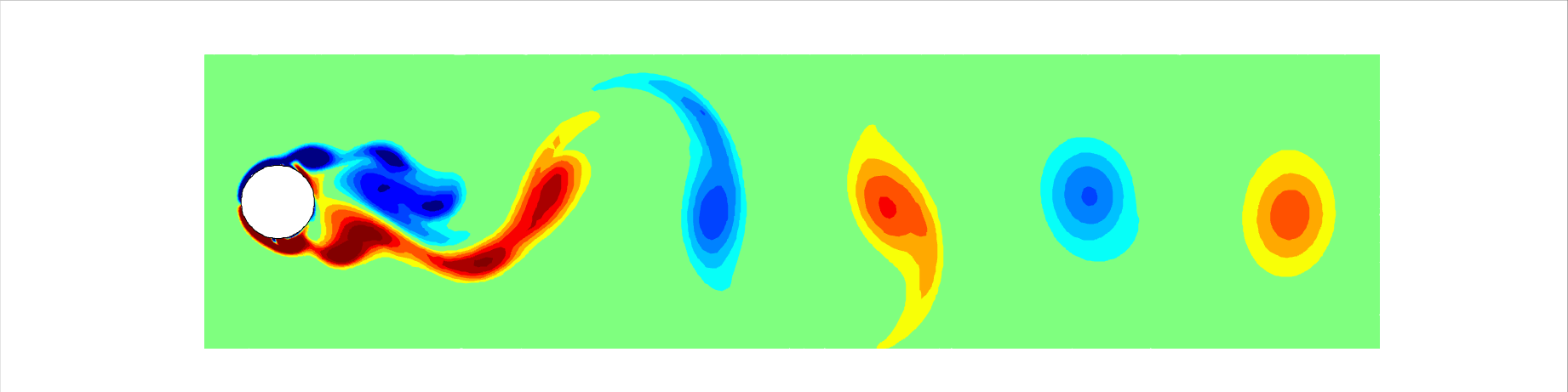}         &
        \rotatebox[origin=c]{90}{$t=81.9$}
        \\[1em]
        
  \includegraphics[trim=8cm 2cm 8cm 2cm,clip,width=0.42\textwidth]{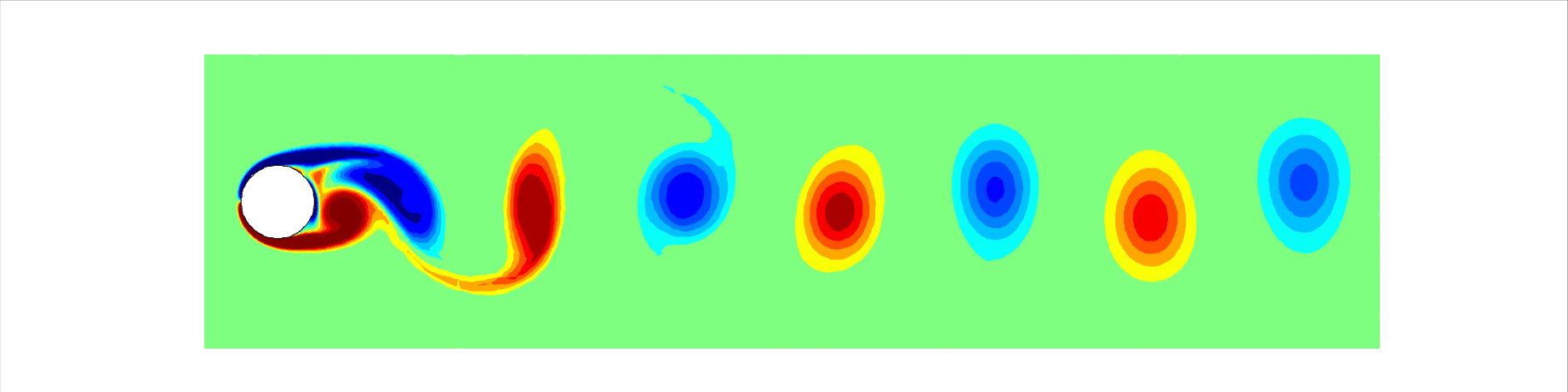} &
        \rotatebox[origin=c]{90}{$t=24.1$}&        
        \includegraphics[trim=8cm 2cm 8cm 2cm,clip,width=0.42\textwidth]{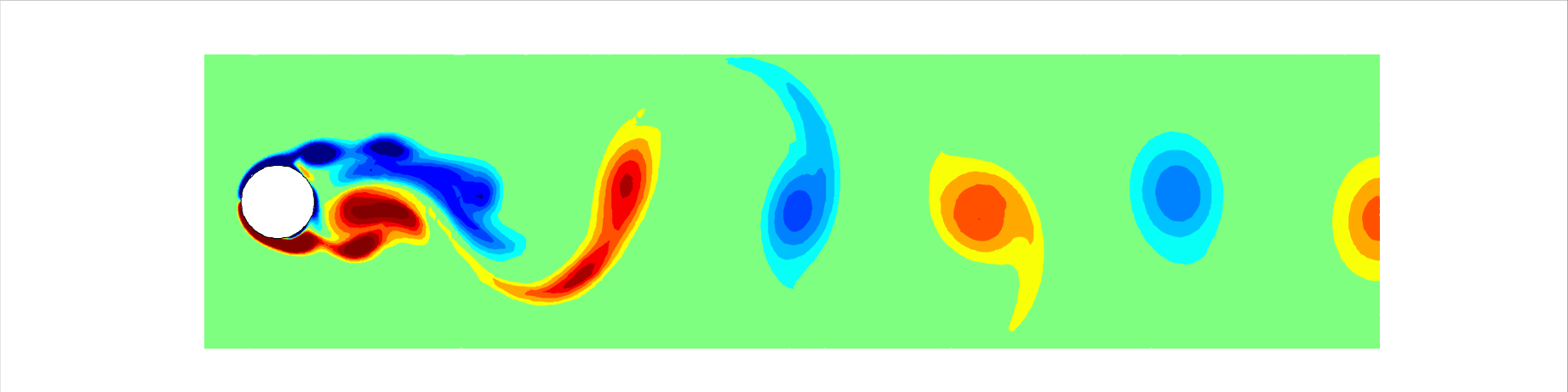}         &
        \rotatebox[origin=c]{90}{$t=83.2$}
        \\[1em]
        
  \includegraphics[trim=8cm 2cm 8cm 2cm,clip,width=0.42\textwidth]{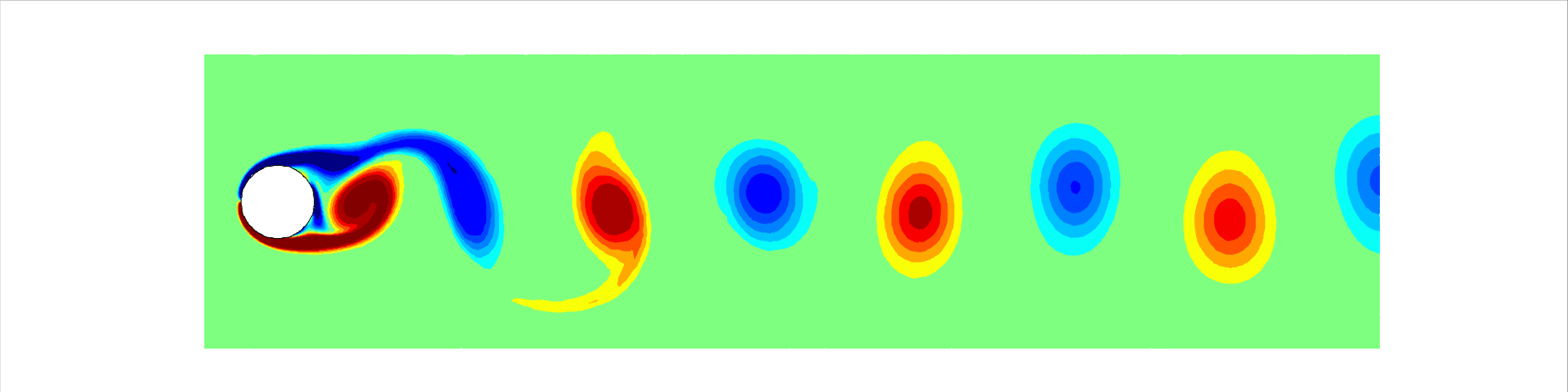} &
        \rotatebox[origin=c]{90}{$t=25.3$}&        
        \includegraphics[trim=8cm 2cm 8cm 2cm,clip,width=0.42\textwidth]{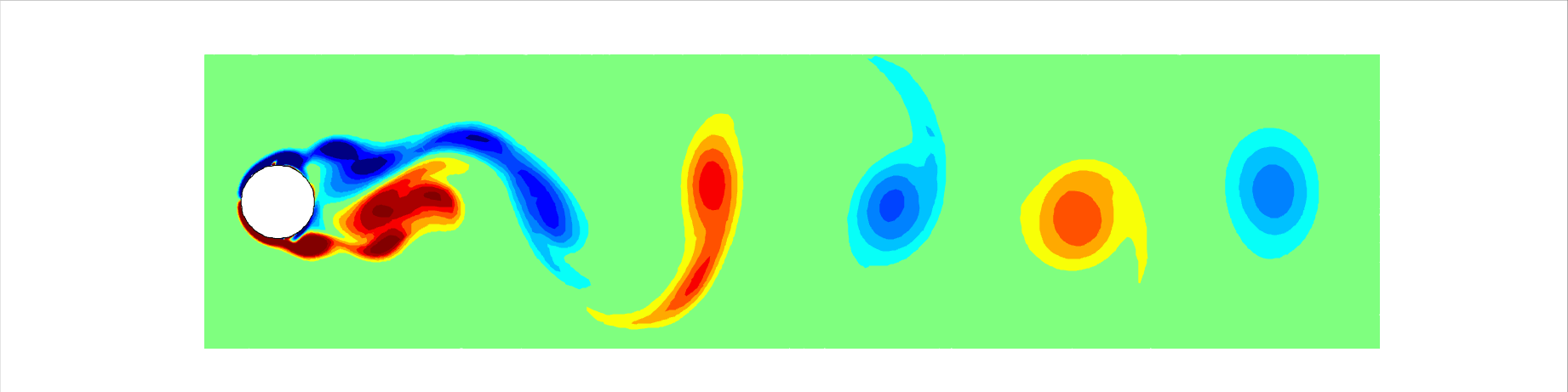}         &
        \rotatebox[origin=c]{90}{$t=84.6$}
        \\[1em]
        
 \includegraphics[trim=8cm 2cm 8cm 2cm,clip,width=0.42\textwidth]{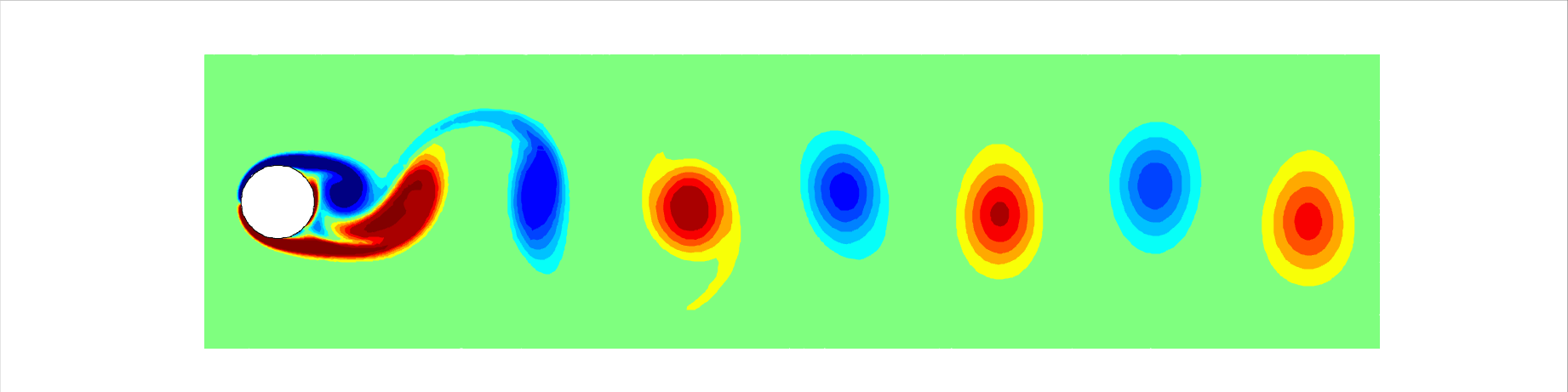} &
        \rotatebox[origin=c]{90}{$t=26.5$}&        
        \includegraphics[trim=8cm 2cm 8cm 2cm,clip,width=0.42\textwidth]{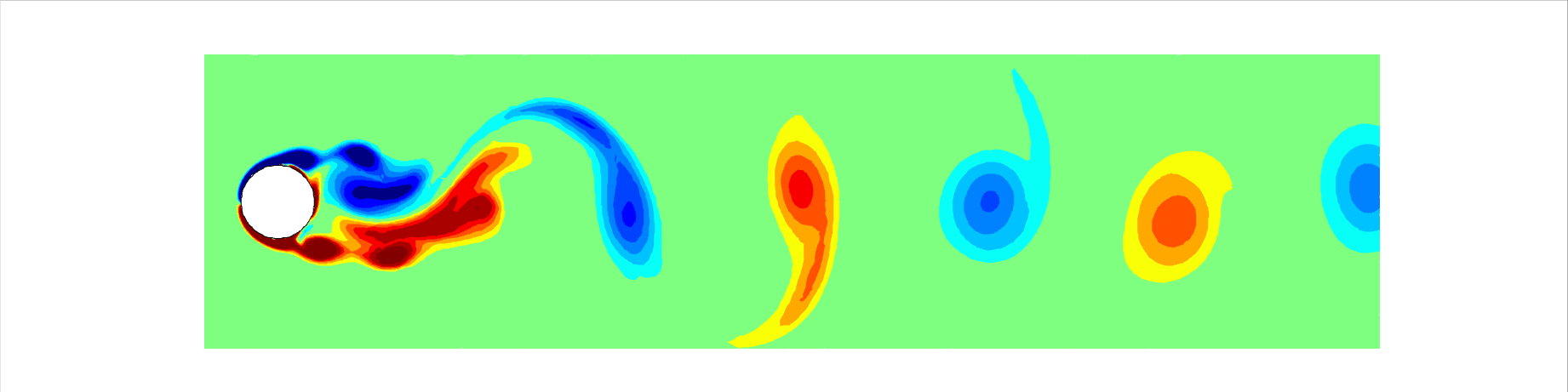}         &
        \rotatebox[origin=c]{90}{$t=86.0$}
        \\[1em]

    \end{tabular}
    \caption{Type-1 control: selected snapshots of the vorticity field near the cylinder at $Re = 300$. Left: uncontrolled baseline case; Right: FML-DRL controlled case.}
    \label{fig:flow_vis_ppo}
\end{figure}

\begin{figure}[ht!]
    \centering
    
    \begin{tabular}{ m{0.42\textwidth} @{}c m{0.42\textwidth} @{}c@{} }
        
        \includegraphics[trim=8cm 2cm 8cm 2cm,clip,width=0.42\textwidth]{./figures/type1/referencevorticity_000219.png} &
        \rotatebox[origin=c]{90}{$t=21.8$}&        
        \includegraphics[trim=8cm 2cm 8cm 2cm,clip,width=0.42\textwidth]{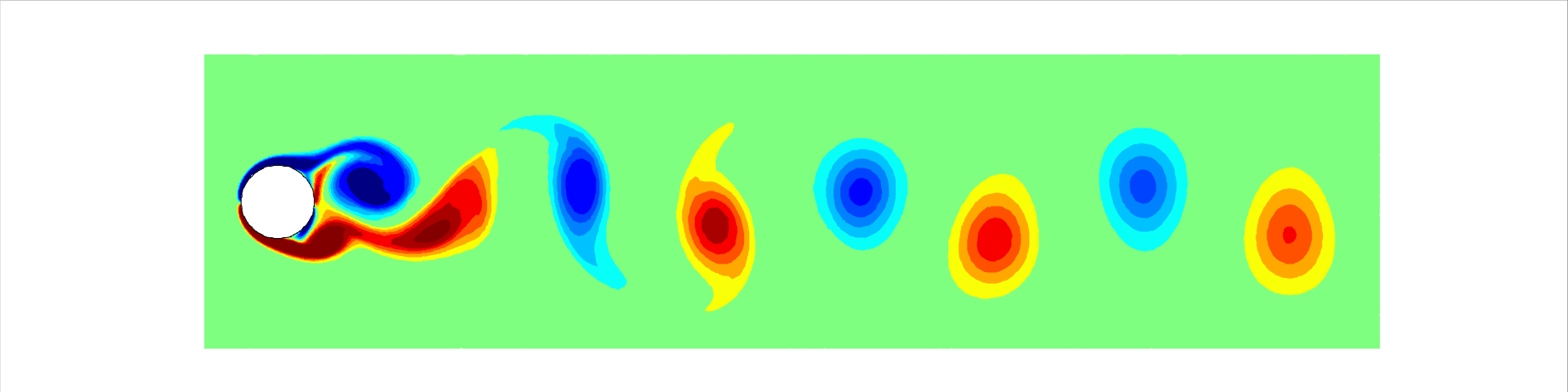}         &
        \rotatebox[origin=c]{90}{$t=81.4$}
        \\[1em]
        
  \includegraphics[trim=8cm 2cm 8cm 2cm,clip,width=0.42\textwidth]{./figures/type1/referencevorticity_000231.png} &
        \rotatebox[origin=c]{90}{$t=23.0$}&        
        \includegraphics[trim=8cm 2cm 8cm 2cm,clip,width=0.42\textwidth]{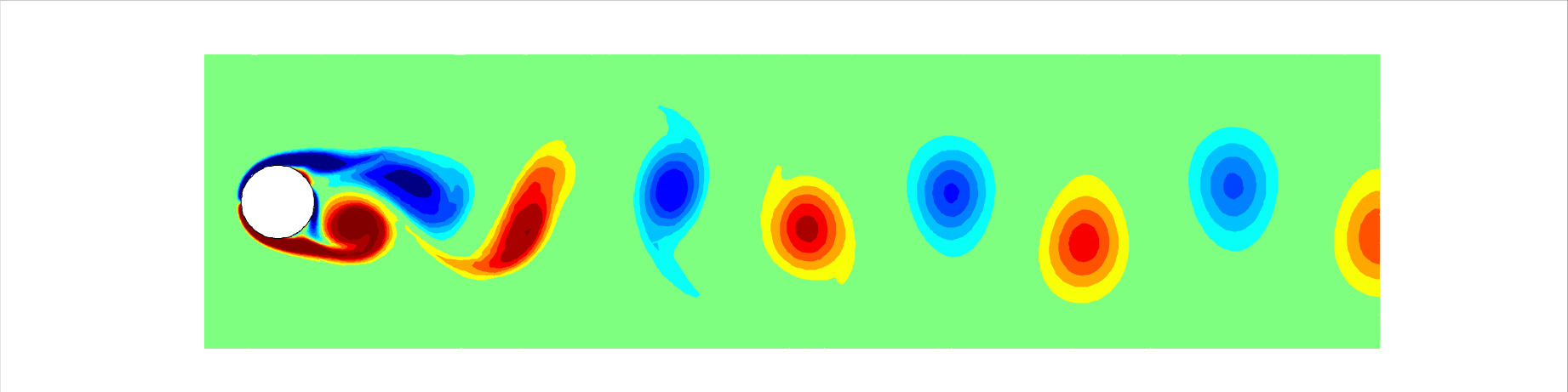}         &
        \rotatebox[origin=c]{90}{$t=82.8$}
        \\[1em]
        
  \includegraphics[trim=8cm 2cm 8cm 2cm,clip,width=0.42\textwidth]{./figures/type1/referencevorticity_000242.png} &
        \rotatebox[origin=c]{90}{$t=24.1$}&        
        \includegraphics[trim=8cm 2cm 8cm 2cm,clip,width=0.42\textwidth]{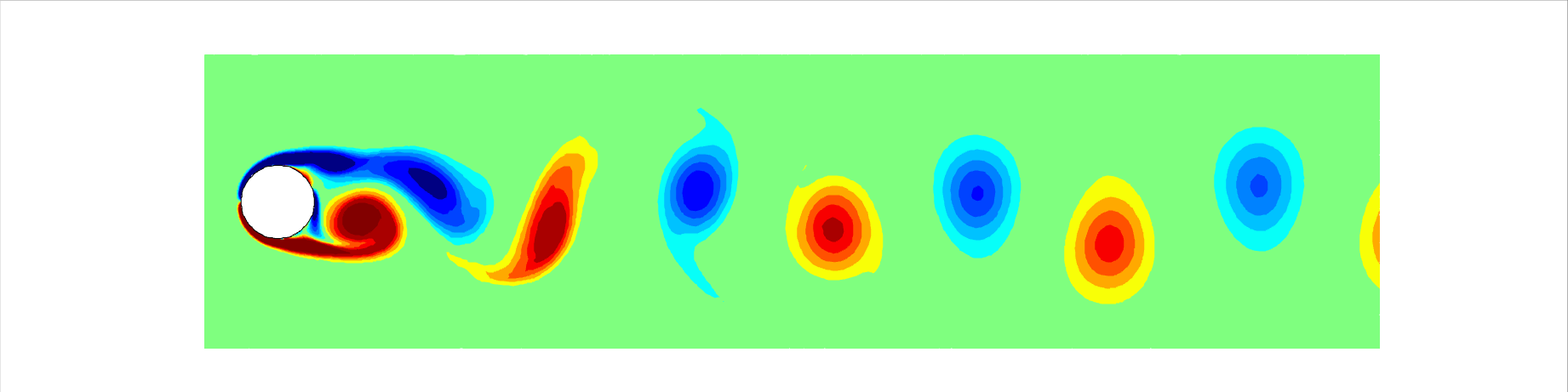}         &
        \rotatebox[origin=c]{90}{$t=83.2$}
        \\[1em]
        
  \includegraphics[trim=8cm 2cm 8cm 2cm,clip,width=0.42\textwidth]{./figures/type1/referencevorticity_000254.png} &
        \rotatebox[origin=c]{90}{$t=25.3$}&        
        \includegraphics[trim=8cm 2cm 8cm 2cm,clip,width=0.42\textwidth]{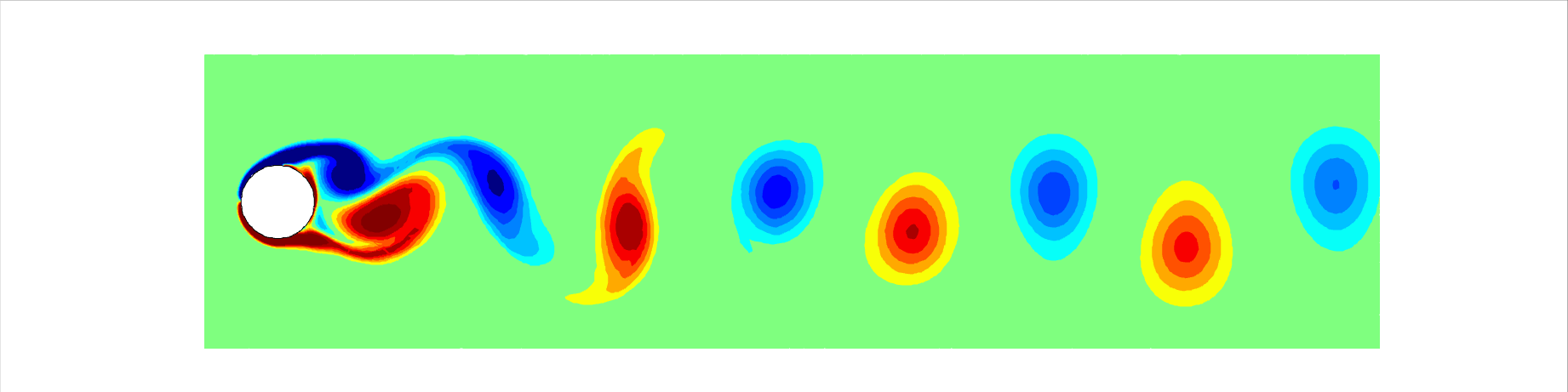}         &
        \rotatebox[origin=c]{90}{$t=84.4$}
        \\[1em]
        
 \includegraphics[trim=8cm 2cm 8cm 2cm,clip,width=0.42\textwidth]{./figures/type1/referencevorticity_000266.png} &
        \rotatebox[origin=c]{90}{$t=26.5$}&        
        \includegraphics[trim=8cm 2cm 8cm 2cm,clip,width=0.42\textwidth]{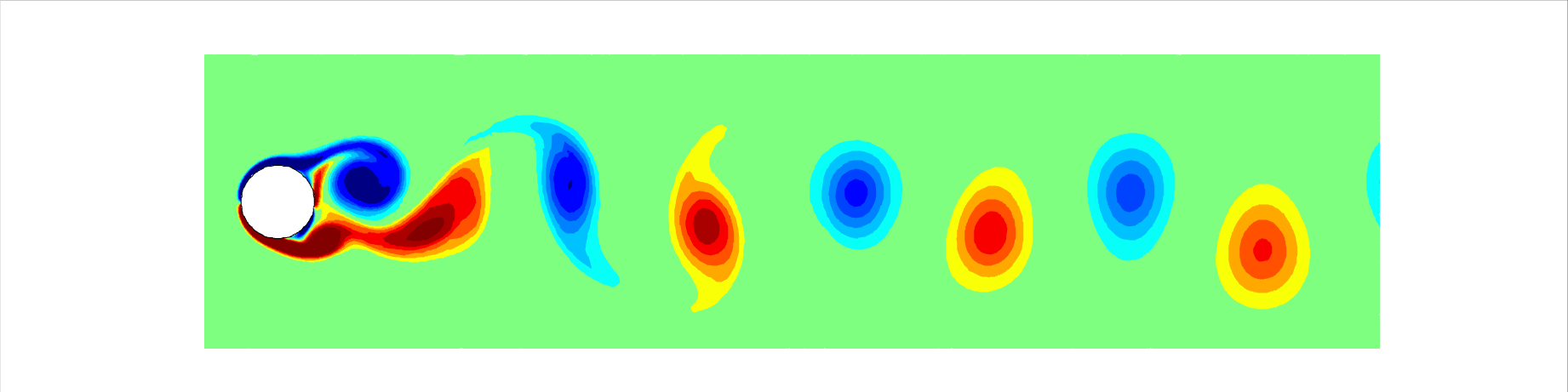}         &
        \rotatebox[origin=c]{90}{$t=85.6$}
        \\[1em]        
    \end{tabular}
     \caption{Type-1 control: selected snapshots of the vorticity field near the cylinder at $Re = 300$. Left: uncontrolled baseline case; Right: FML-MPC controlled case.}
    \label{fig:flow_vis_mpc}
\end{figure}

Fig.~\ref{fig:closed_loop_case1} shows the performance of both DRL-PPO and MPC controllers, compared to the uncontrolled baseline case (which is the standard flow past cylinder simulation). The top of Fig.~\ref{fig:closed_loop_case1} shows the drag coefficient over a time horizon of $T=200$. Here the dashed lines are the instantaneous solutions, whereas the solid lines are the time averaged solutions over a moving window of $T=5.0$. The uncontrolled baseline case naturally produces periodic drag, corresponding to the periodic vortex shedding, with its time averaged value approximately a constant. For the active control case, we observe that both the DRL-PPO controller and the MPC controller produce significant reduction in drag, at approximately 20\% compared to the baseline case. 
The instantaneous control signal, which updates the boundary condition \eqref{eq:jet_control} at the cylinder surface, are shown in the middle of Fig.~\ref{fig:closed_loop_case1}. We observe that DRL-PPO and MPC produce rather different control signals, even though the end drag reduction effects are similar. The lift coefficients, which are not the primary goal of the flow control, are shown in the bottom of the figure, for the sake of completeness.
For the same reason, we present the vorticity field evolution in Fig.~\ref{fig:flow_vis_ppo} for the DRL-PPO control and in Fig.~\ref{fig:flow_vis_mpc} for the MPC control, along with those for the uncontrolled baseline case for comparison. 

\subsection{Type-II Active Flow Control}

For Type-II control, we start the flow from a fully developed stable state at arbitrarily chosen $Re$ number and then apply the FML-2 model (without using the $Re$ information) and start the MPC controller (Section \ref{sec:MPC}). Note that the FML-2 model was trained with unknown $Re$ in range $Re\in [100,500]$.

\begin{figure}[ht!]
    \centering
    \begin{subfigure}[b]{0.475\textwidth}
        \centering
        \includegraphics[width=\textwidth]{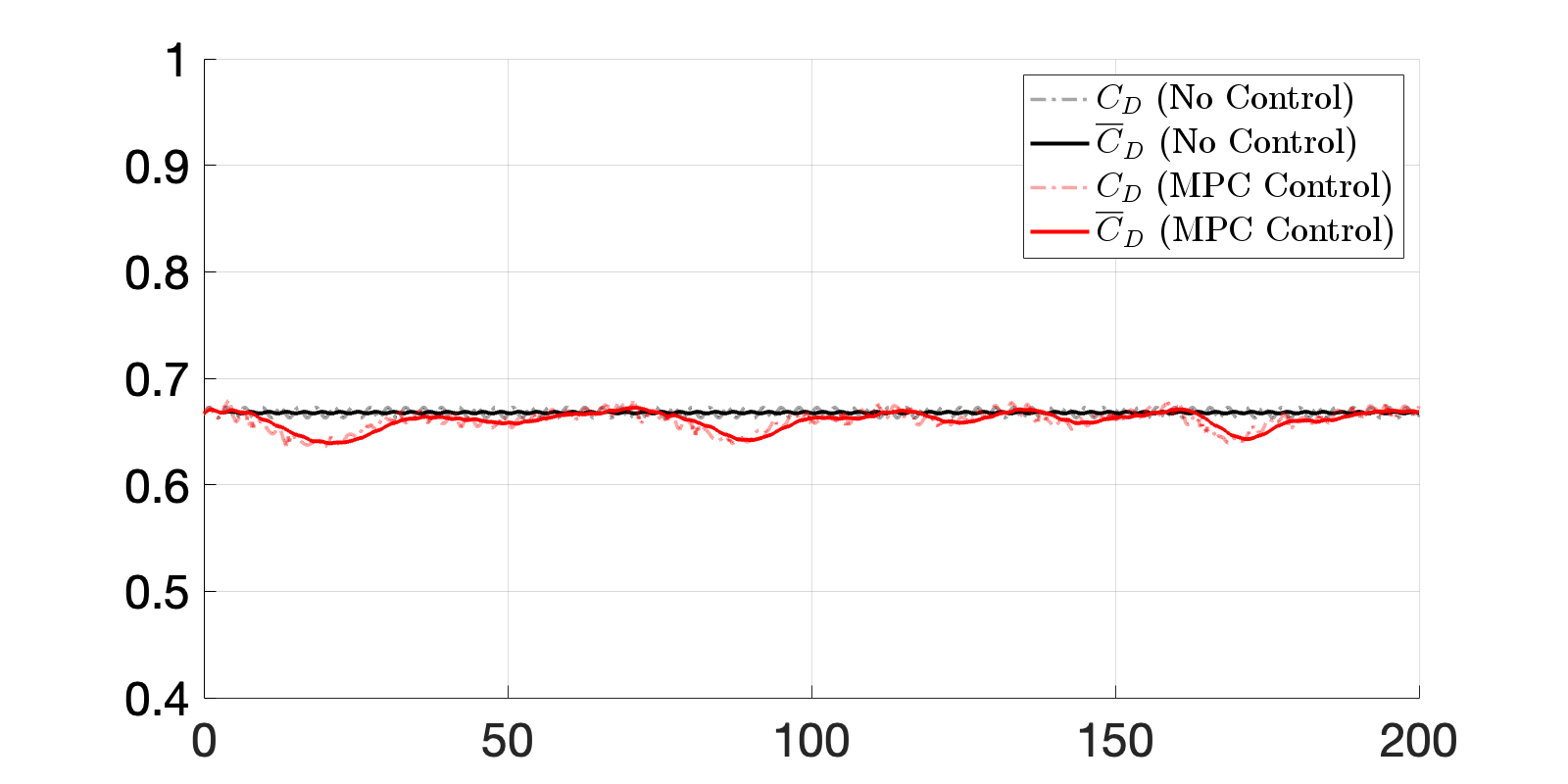}
        \caption{$Re=100$}
    \end{subfigure}
    \hfill
    \begin{subfigure}[b]{0.475\textwidth}
        \centering
        \includegraphics[width=\textwidth]{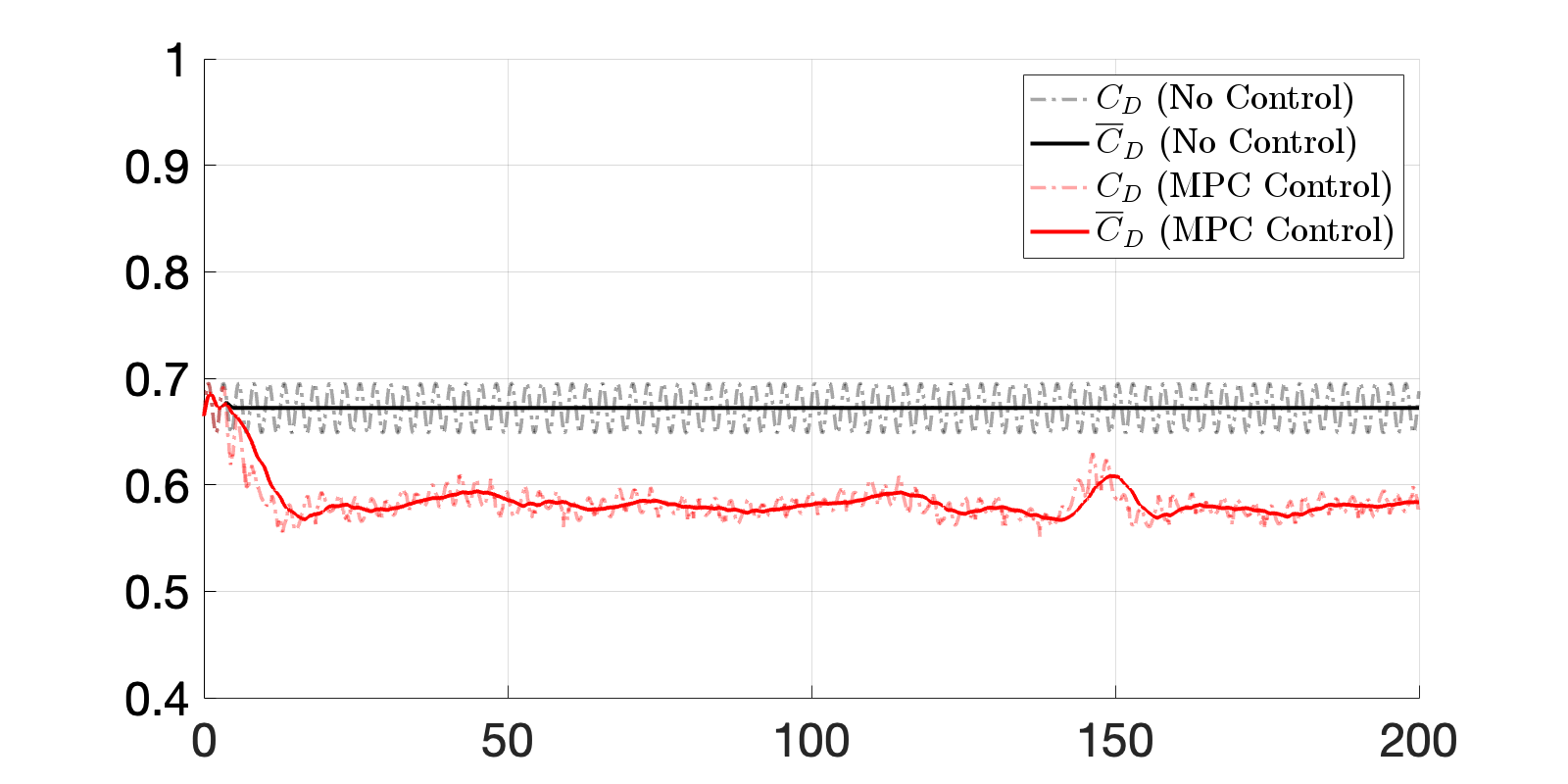}
        \caption{$Re=200$}
    \end{subfigure}

    \begin{subfigure}[b]{0.475\textwidth}
        \centering
        \includegraphics[width=\textwidth]{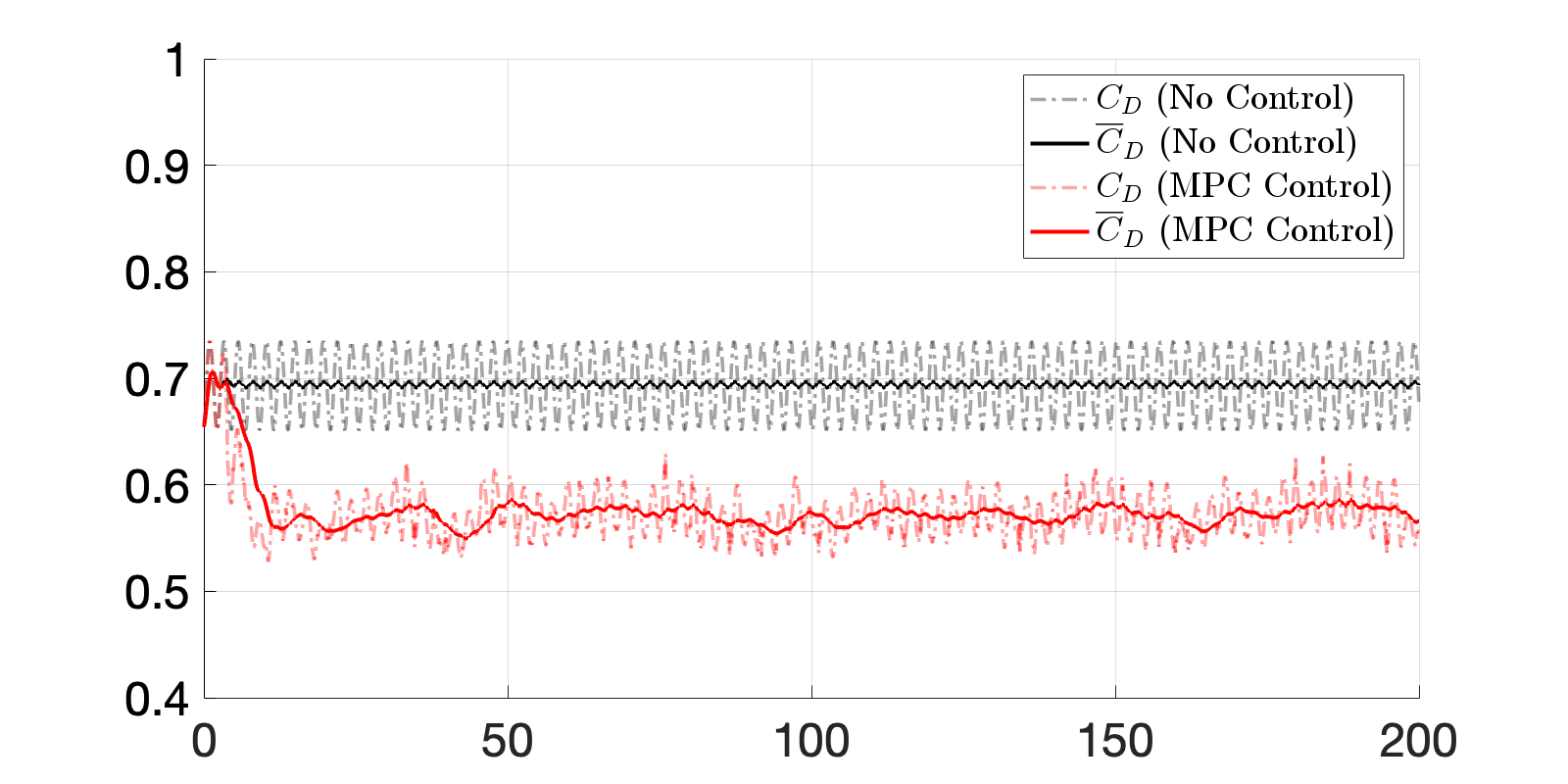}
        \caption{$Re=300$}
    \end{subfigure}
    \hfill
    \begin{subfigure}[b]{0.475\textwidth}
        \centering
        \includegraphics[width=\textwidth]{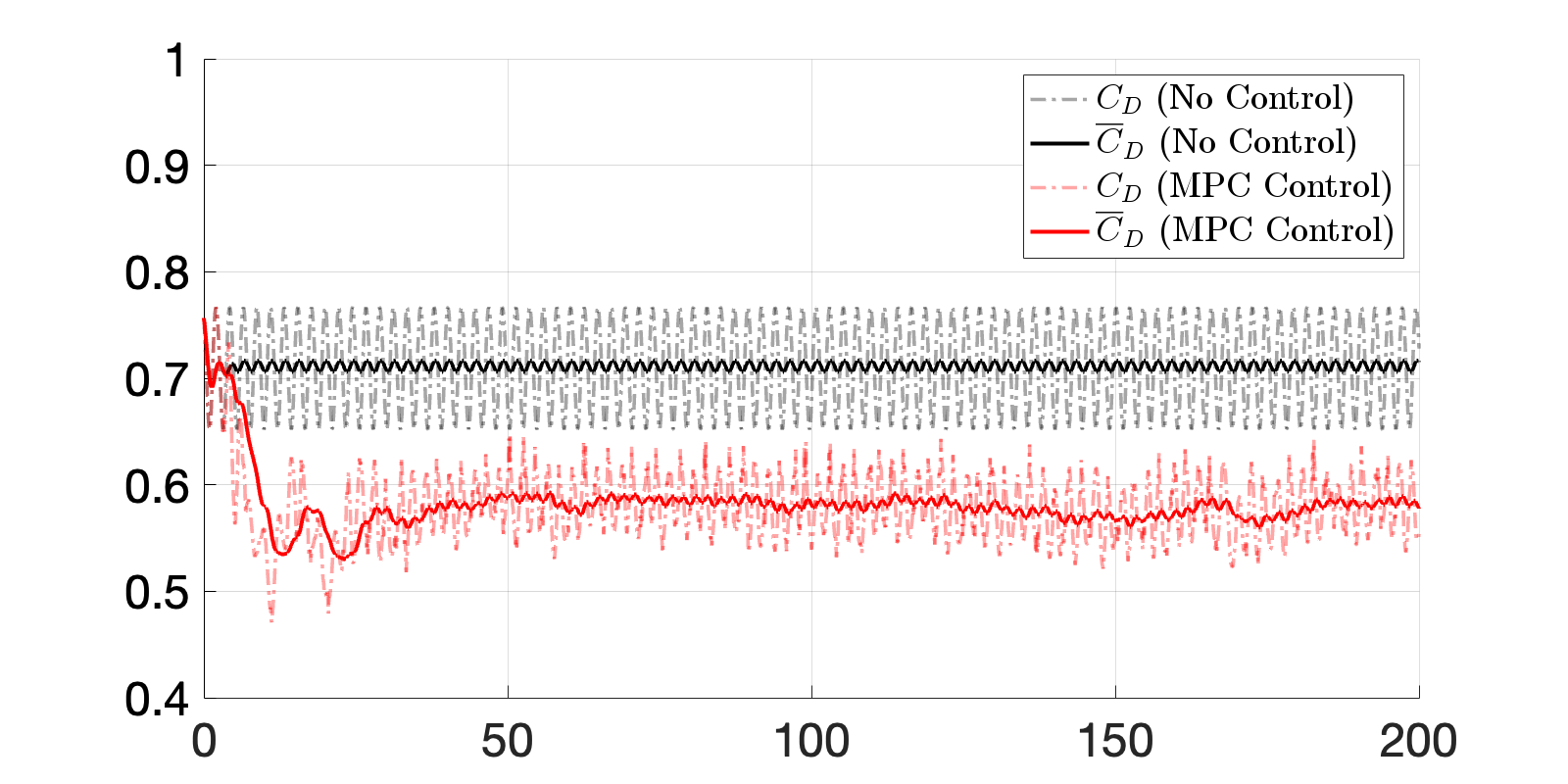}
        \caption{$Re=400$}
    \end{subfigure}

    \begin{subfigure}[b]{0.475\textwidth}
        \centering
        \includegraphics[width=\textwidth]{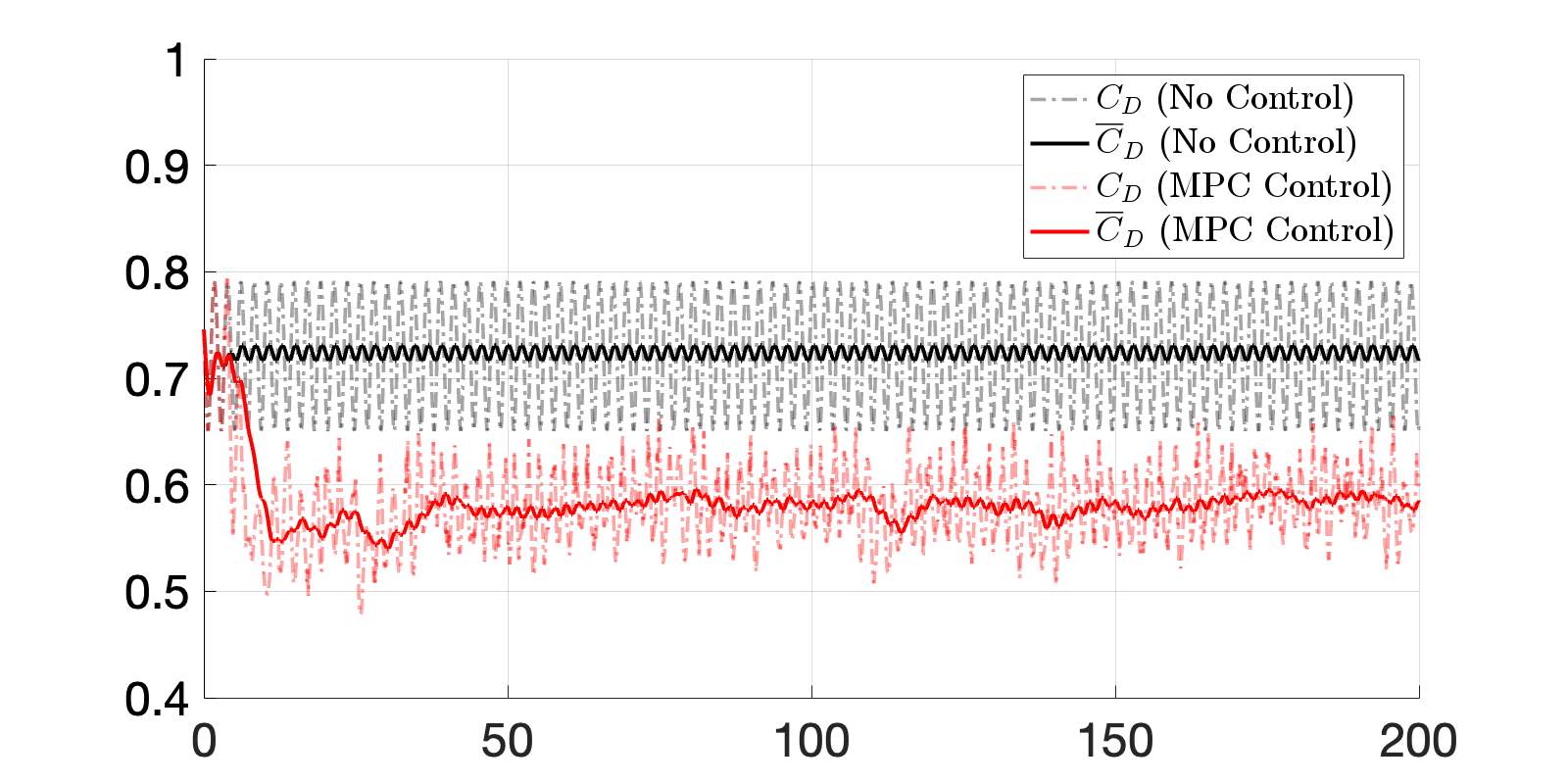}
        \caption{$Re=500$}
    \end{subfigure}
    \hfill
    \begin{subfigure}[b]{0.475\textwidth}
        \centering
        \includegraphics[width=\textwidth]{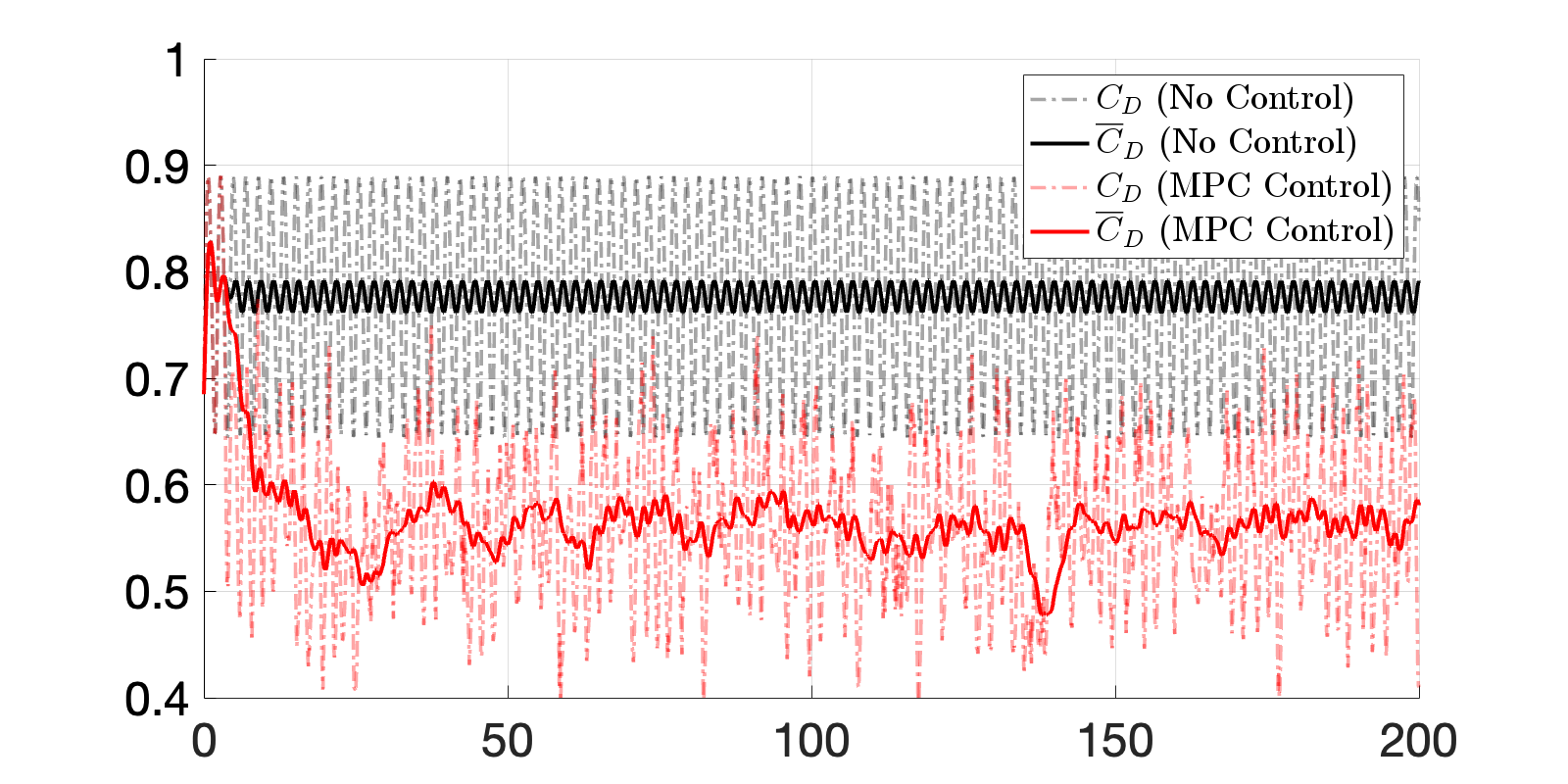}
        \caption{$Re=1000$}
    \end{subfigure}
    \caption{Type-II MPC optimal control using FML-2 model at unknown $Re$ numbers, for 6 test cases of $Re=100, 200, 300, 400, 500, 1000$. Dashed lines are instantaneous drag coefficients; Solid lines are averaged values over moving window of size $5.0$. Note that $Re$ number is a hidden parameter and does not enter the FML-2 model.}
    \label{fig:generalization}
\end{figure}

In Fig.~\ref{fig:generalization} we present 6 different cases for the MPC control. The 6 cases are $Re=100$, $200$, $300$, $400$, $500$, and $1,000$. We remark that the $Re=1,000$ case is completely outside the training regime for the FML-2 model, whereas $Re=100$ and $Re=500$ are borderline cases. The results demonstrate that the FML-2 MPC controller is able to produce noticeable drag reduction in all cases. The higher the $Re$ number, the more reduction the MPC control achieves, which is not surprising. Again, we emphasize that the FML-2 model is completely unaware of the $Re$ number in each case.

\section{Conclusion}

This work presents a framework for optimal control based on learning the dynamics of quantities of interest (QoIs) using flow map learning. By modeling QoI dynamics directly, the approach avoids repeated solution of the full-order governing equations during control optimization and enables fast evaluation of control strategies. 

Using active flow control (drag reduction) for incompressible flow past a cylinder as a representative test case, we demonstrate that accurate models of QoI dynamics can be learned directly and used for control strategies such as deep reinforcement learning (DRL) and model predictive control (MPC). Without needing the full-scale flow solver, the proposed method enables on-the-fly active flow control and offers the potential for real-time optimal control for complex systems.

\bibliographystyle{siamplain}
\bibliography{reference}

\end{document}